\documentclass{amsproc}%
\usepackage{amsfonts}%
\usepackage{amsmath}%
\usepackage{amssymb}%
\usepackage{graphicx}
\theoremstyle{plain}

\newtheorem{corollary}{Corollary}

\newtheorem{definition}{Definition}
\newtheorem{example}{Example}

\newtheorem{lemma}{Lemma}

\newtheorem{proposition}{Proposition}

\newtheorem{theorem}{Theorem}
\numberwithin{equation}{section}
\begin{document}
\title{The Deformation Complex for DG Hopf Algebras }
\author{Ronald N. Umble}
\address{Department of Mathematics\\
Millersville University\\
Millersville, PA 17551\\
(717) 872-3708}
\email{ron.umble@millersville.edu}
\urladdr{http://www.millersville.edu/\symbol{126}rumble}
\thanks{This research funded in part by a Millersville University Faculty Reasearch Grant.}
\date{October 1, 1995}
\subjclass{16W30, 13D10, 16E45, 55P62}
\keywords{Hopf Algebra, $A\left(  n\right)  $-algebra, deformation}

\begin{abstract}
Let $H$ be a DG Hopf algebra over a field $\mathbf{k}$. This paper gives an
explicit construction of a triple cochain complex that defines the
Hochschild-Cartier cohomology of $H$. A certain truncation of this complex is
the appropriate setting for deforming $H$ as an $H\left(  q\right)
$-structure. The direct limit of all such truncations is the appropriate
setting for deforming $H$ as a strongly homotopy associative structure. Sign
complications are systematically controlled. The connection between rational
perturbation theory and the deformation theory of certain free commutative
differential graded algebras is clarified.

\end{abstract}
\maketitle

\section{Introduction\bigskip\ }

The purpose of this paper is two-fold: (1) to give an explicit construction of
the deformation complex for differential graded Hopf algebras and (2) to
relate the rational perturbation theory of Felix \cite{Felix} and
Halperin-Stasheff \cite{Halperin} to the deformation theory of certain free
commutative differential graded algebras. The untruncated deformation complex
constructed here directs the deformation of a differential graded Hopf algebra
$H$ as an $H\left(  \infty\right)  $-structure; appropriate truncations direct
the deformation of $H$ as an $H(q)$-structure. The special case $q=3$ is
applied by Lazarev and Movshev in their paper \textit{Deformations of the de
Rham Algebra }\cite{Lazarev}, which follows as a sequel.

In \cite{Gersten}, Gerstenhaber and Schack showed how to deform a
biassociative Hopf algebra $H$ over a field $\mathbf{k}$ relative to its
algebraic cohomology. Following their cues, we define the algebraic cohomology
of a connected biassociative differential graded Hopf algebra $H$ and give a
brief exposition of the related deformation theory.

This exposition minimizes the sign complications that arise in a graded theory
by adopting two strategies: (1) we work at the (coordinate free) operator
level and (2) we base our constructions on $H$-free resolutions with
differentials of internal degree zero. Thus, elements of $H$ never ''move
past'' graded cochains and graded cochains are free to ''move past'' the
resolution differentials without complicating signs. While the first strategy
is evident in Gerstenhaber and Schack's exposition \cite{Gersten}, the second
was used by Burghelea and Poirrier \cite{Burghelea} to define the Hochschild
and cyclic cohomologies of free commutative associative differential graded
algebras in characteristic zero. The recent work of Penkava and Schwarz
\cite{Penkava} demonstrates that careful attention to signs can be critical.

This paper is organized as follows: Section 2 establishes the necessary
preliminaries and section 3 reviews the ''classical'' (co)bar resolution
\cite{MacLane} of a graded (co)algebra and its extension to a differential
graded (co)algebra. These resolutions, which are \textit{not} meant to model
chains on some contractible space, have differentials of internal degree zero
and avoid the dimension shifts of Adams \cite{Adams} and Eilenberg and Mac
Lane \cite{Eilenberg}. Section 4 dualizes and generalizes the notion of a
differential graded bimodule over a differential graded algebra, which is
implicit in \cite{Markl}, to analogous structures over differential graded
coalgebras and Hopf algebras.

In section 5 we define the Hochschild cohomology of a connected associative
differential graded algebra (\textit{d.g.a.) }$A$ with coefficients in a
differential graded $A$-bimodule $M$. The deformation complex for $A$ is
obtained by setting $M=A$ and appropriately truncating the Hochschild cochain
complex$.$ A construction of this cohomology was given by Markl in
\cite{Markl} but with two significant differences: (1) our underlying bar
resolution does not use the Eilenberg-Mac Lane dimension shift and (2) we
transfer the theory from the level of $A$-bimodules to the level of
$\mathbf{k}$-modules at which the cohomology and deformation theory are
clearly linked.

We also define the Harrison cohomology of a commutative \textit{d.g.a.
(c.d.g.a) }with coefficients in a symmetric\textit{\ d.g. }$A$-bimodule $M$.
The Harrison cohomology of \textit{free c.d.g.a.'s with trivial coefficients}
was defined earlier by Burghelea and Poirrier \cite{Burghelea}. We show how to
interpret the rational perturbation theory of ''bigraded models''
\cite{Halperin}, which are certain free \textit{c.d.g.a}.'s over
$\mathbb{Q}\mathbf{,}$ in terms of the ''appropriately truncated'' Harrison
cohomology of the model with coefficients in itself. We observe that for free
\textit{c.d.g.a.'s, }all flexibility lies in the direction of the
differential. The Lie algebra analogs of these constructions recently appeared
in \cite{Lada}.

Next we dualize and obtain the Cartier\textit{\ }cohomology of a connected
coassociative differential graded coalgebra (\textit{d.g.c.) }$C$ with
coefficients in a differential graded $C$-bicomodule $N;$ the deformation
complex for $C$ is obtained by setting $N=C$ and appropriately truncating the
Cartier cochain complex. Finally, we join these dual theories and obtain the
Hochschild-Cartier\textit{\ }cohomology of a connected biassociative
differential graded Hopf algebra (\textit{d.g.h.a.) }$H$\textit{;} the
deformation complex for $H$ is an appropriate truncation of the triple cochain
complex for this cohomology. Section 6 concludes the discussion with a brief
exposition of the deformation theory for \textit{d.g.h.a.'s}.\ 

\section{Notation and Preliminaries\ }

Let $R$ be a commutative ring with identity $1_{R}$ and let $M$ be a
(non-negatively) graded $R$-module. $M$ is \textit{connected }if $M^{0}\approx
R.$ Unless indicated otherwise, all tensor products will be defined over $R.$
Let $\{M_{i}\}$ be a sequence of graded $R$-modules; the subspace of
$\otimes_{i}M_{i}$ consisting of all elements homogeneous in degree $p$ is
deonoted by $(\otimes_{i}M_{i})^{p}.$ Let $M^{\otimes n}=M\otimes\cdots\otimes
M$ with $n>0$ factors and define $M^{\otimes0}=R.$ Let $T_{p,n}M=(M^{\otimes
n})^{p},$ then $TM=\sum_{p,n\geq0}T_{p,n}M$ is a bigraded space; a
bihomogeneous element $x\in TM$ has bidegree $(p,n)$ and is said to have
\textit{internal degree p and external degree n.} The symbol $\mid x\mid$
denotes the \textit{internal} degree of $x$\textit{. }

A map $f:M\rightarrow M$ has \textit{degree} $p$ if $\mid f(x)\mid=\mid
x\mid+p$ for each homogeneous $x\in M,$ in which case we write $\mid f\mid=p.$
The identity map $1:M\rightarrow M$ and the canonical isomorphisms
$i_{1}:M\rightarrow R\otimes M,$ $j_{1}:R\otimes M\rightarrow M,$
$i_{2}:M\rightarrow M\otimes R,$ and $j_{2}:M\otimes R\rightarrow M$ are maps
of degree zero. Another such map is the \textit{permutation operator}
$\sigma:M^{\otimes n}\rightarrow M^{\otimes n},$ defined by $\sigma
(x_{1}\otimes\cdots\otimes x_{n})=$ $\pm\ x_{\sigma^{-1}(1)}\otimes
\cdots\otimes x_{\sigma^{-1}(n)},$ where $\sigma\in S_{n}$ and the sign is
given by the standard sign commutation rule with respect to \textit{internal}
degree\textit{: whenever two symbols} $u$ \textit{and} $v$ \textit{with
internal degrees are interchanged, affix the sign }$(-1)^{\mid u\mid\mid
v\mid}$\textit{\ }(see \cite{MacLane}, p. 164). If $x,y\in M$ and
$f,g:M\rightarrow M$, the sign commutation rule gives: $(f\otimes g)(x\otimes
y)=(-1)^{\mid g\mid\mid x\mid}f(x)\otimes g(y).$ An $R$-module map
$d:M\rightarrow M$ of degree $\pm1$ such that $d\circ d=0$ is called
a\textit{\ differential} on $M;$ the pair $(M,d)$ is a \textit{differential
graded (d.g.) }$R$-module.

Let $A$ be a graded symmetric $R$-module. A \textit{multiplication on }$A$ is
an $R$-module map $\mu:A\otimes A\rightarrow A$ of degree $0$; the pair
$(A,\mu)$ is a \textit{graded $R$-algebra}. An $R$-algebra $(A,\mu)$ is
\textit{associative }if $\mu\circ(\mu\otimes1)=\mu\circ(1\otimes\mu);$ it is
\textit{commutative }if $\mu=\mu\circ(1,2).$ It is \textit{unital} if there
exists an $R$-algebra map $\eta:R\rightarrow A$ of degree $0$ such that
$\mu\circ(\eta\otimes1)\circ i_{1}=\mu\circ(1\otimes\eta)\circ i_{2}=1,$ in
which case $\eta$ is called the \textit{unit.} The element $1_{A}=\eta(1_{R})$
acts as a two-sided identity for $\mu$. When $A$ is connected, the unique
algebra isomorphism $R\rightarrow A^{0}$ is a canonical unit$.$ A
\textit{derivation }of\textit{\ }$(A,\mu)$ is an $R$-module map $\theta
:A\rightarrow A $ such that $\theta\circ\mu=\mu\circ(\theta\otimes
1+1\otimes\theta). $ If a differential $d$ on $A$ is a derivation of $(A,\mu
)$, then $d$ is an \textit{algebra differential} and the triple $(A,d,\mu)$ is
a \textit{differential graded }$R$\text{\textit{-algebra (d.g.a.).}}

Let $n\in\mathbf{N}\cup\{\infty\}.$ An $A\left(  n\right)  $-\textit{algebra}
is defined to be a graded $R$-module $A$ together with maps $\{\mu^{(\ell)}\in
Hom_{R}^{2-\ell}(A^{\otimes\ell},A)\}_{1\leq\ell\leq n}$ such that for each
$\ell\leq n,$
\[
\sum_{0\leq i<j;\text{ }j+k=\ell+1}(-1)^{i+ik+\ell k+k}\mu^{(j)}%
\circ(1^{\otimes i}\otimes\mu^{(k)}\otimes1^{\otimes(j-i-1)})=0.
\]
The signs here agree with those in \cite{Stasheff}; we use upper indices and
reserve the lower for indexing coefficients in a deformation. An
$A(n)$-algebra is \textit{strict} if $\mu^{(n)}=0.$ Every \textit{d.g.a.
}$(A,d,\mu)$ is a strict $A(n)$-algebra for all $n\geq3$ via $\mu
^{(1)}=d,\ \mu^{(2)}=\mu,\ $and $\mu^{(i)}=0$ for $3\leq i\leq n.$

Let $C$ be a graded symmetric $R$-module. A \textit{comultiplication}
\textit{on }$C$ is an $R$-module map $\Delta$$:C\rightarrow C\otimes C$ of
degree $0$; the pair $(C,\Delta)$ is a \textit{graded $R$-coalgebra}. An
$R$-coalgebra $(C,\Delta)$ is \textit{coassociative} if $(\Delta\otimes
1)\circ\Delta=(1\otimes\Delta)\circ\Delta;$ it is \textit{cocommutative} if
$\Delta=(1,2)\circ\Delta$. It is \textit{counital} if there exists an
$R$-coalgebra map $\varepsilon:C\rightarrow R$ of degree $0$ such that
$j_{1}\circ(\varepsilon\otimes1)\circ\Delta=j_{2}\circ(1\otimes\varepsilon
)\circ\Delta=1,$ in which case $\varepsilon$ is called the \textit{counit}.
When $A$ is connected, the unique coalgebra isomorphism $C^{0}\rightarrow R$
extended to the zero map in positive degrees is a canonical counit$.$ A
\textit{coderivation of }$(C,\Delta)$ is an $R$-module map $\omega
:C\rightarrow C$ such that $\Delta\circ\omega=(\omega\otimes1+1\otimes
\omega)\circ\Delta.$ If a differential $d$ on $C$ is a coderivation of
$(C,\Delta),$ then $d$ is a \textit{coalgebra differential} and the triple
$(C,d,\Delta)$ is a \textit{differential graded }$R$\textit{-coalgebra
(d.g.c.).}

Let $m\in\mathbf{N}\cup\{\infty\}.$ An $C\left(  m\right)  $%
-\textit{coalgebra} is defined to be a graded $R$-module $C$ together with
maps $\{\Delta^{(\ell)}\in Hom_{R}^{2-\ell}(C,C^{\otimes\ell})\}_{1\leq
\ell\leq m}$ such that for each $\ell\leq m,$
\[
\sum_{0\leq i<j;\text{ }j+k=\ell+1}(-1)^{i+ik+\ell k+k}(1^{\otimes i}%
\otimes\Delta^{(k)}\otimes1^{\otimes(j-i-1)})\circ\Delta^{(j)}=0.
\]
An $C(m)$-coalgebra is \textit{strict if }$\Delta^{(m)}=0.$ Every
\textit{d.g.c.} $(C,d,\Delta)$ is a strict $C(m)$-coalgebra for all $m\geq3.$

Let $H$ be a graded symmetric $R$-module, and suppose that $H$ is equipped
with a multiplication $\mu,$ a unit $\eta,$ a comultiplication $\Delta,$ and a
counit $\varepsilon$ such that $\eta$ and $\varepsilon$ are $R$-bialgebra maps
and $\Delta\circ\mu=(\mu\otimes\mu)\circ(2,3)\circ(\Delta\otimes\Delta);$ then
$(H,\mu,\eta,\Delta,\varepsilon)$ is a \textit{graded }$R$-\textit{bialgebra}.
This latter condition is equivalent to requiring that $\mu$ and $\Delta$ be,
respectively, coalgebra and algebra maps. An \textit{antipode} for a graded
$R$-bialgebra $H$ is an $R$-antialgebra map $S:H\rightarrow H$ of degree $0$
such that $\mu\circ(S\otimes1)\circ\Delta=\mu\circ(1\otimes S)\circ\Delta
=\eta\circ\varepsilon.$ A graded $R$-bialgebra $H$ is \textit{biassociative}
if it is both associative and coassociative. When $H$ is connected and
biassociative, there is a unique inductively defined antipode $S$ that acts as
the identity in degree $0$ and by $S(x)=-x-\sum x_{(1)}S(x_{(2)})$ in positive
degrees, where $\Delta(x)=\sum x_{(1)}\otimes x_{(2)};$ see \cite{Milnor}$.$ A
graded $R$-bialgebra $(H,\mu,\eta,\Delta,\varepsilon)$ equipped with antipode
$S$ is a \textit{graded }$R$-\textit{Hopf algebra} \textit{(g.h.a.).}
Furthermore, if $(H,d,\mu)$ is a \textit{d.g.a}., $(H,d,\Delta)$ is a
\textit{d.g.c}. and $(H,\mu,\eta,\Delta,\varepsilon,S)$ is a \textit{g.h.a}.,
then $(H,d,\mu,\eta,\Delta,\varepsilon,S)$ is a \textit{differential graded
}$R$\text{\textit{-Hopf algebra (d.g.h.a.). }}

Henceforth, all objects are assumed to be graded; all $R$-modules are assumed
to be connected; all $R$-algebras, $R$-coalgebras, and $R$-Hopf algebras are
assumed to be associative, coassociative, and biassociative, respectively. An
$R$-Hopf algebra will be unambigously denoted by $(H,\mu,\Delta)$.\ \ 

\section{Two-sided Bar and Cobar Resolutions\ }

\subsection{The Bar Resolution}

Let $\mathbf{k}$ be a field and let $(A,\mu)$ be a $\mathbf{k}$-algebra. For
each $m\geq0,$ inductively define $\mathbf{k}$-linear maps $\partial
_{(m)}:A^{\otimes(m+2)}\rightarrow A^{\otimes(m+1)}$ by
\[
\partial_{(0)}=\mu,
\]
and
\[
\partial_{(m)}=\mu\otimes1^{\otimes m}-1\otimes\partial_{(m-1).}
\]
In more familiar form this is
\begin{equation}
\label{bar}\partial_{(m)}=\sum_{i=0}^{m}(-1)^{i}(1^{\otimes i}\otimes
\mu\otimes1^{\otimes(m-i)}),
\end{equation}
but many of the facts we need flow more easily from the inductive form. Let
$\partial=\sum_{m\geq0}\partial_{(m)};$ using induction and the fact that
$\mu$ is associative, it is a simple matter to show that $\partial
\circ\partial=0.$ Hence $\partial$ is a differential with respect to external
degree. The chain complex
\[
A\overset{\partial_{(0)}}{\leftarrow}A\otimes A\overset{\partial_{(1)}%
}{\leftarrow}A\otimes A\otimes A\overset{\partial_{(2)}}{\leftarrow}%
A^{\otimes4}\overset{\partial_{(3)}}{\leftarrow}\cdots
\]
is called the (\textit{classical) two-sided bar resolution of A }%
\cite{MacLane}\textit{.} Furthermore, this resolution is acyclic via
contracting homotopy $s=\sum_{m\geq-1}s_{m}$ where $s_{m}=[(\eta\otimes1)\circ
i_{1}]\otimes1^{\otimes(m+1)}.$ Note that maps $\partial$ and $s$ have degree
zero with respect to the internal grading.

The bar resolution extends to a \textit{d.g.a. }$(A,d,\mu)$ as follows. For
each $m\geq-1,$ inductively define $\mathbf{k}$-linear maps $d_{(m)}%
:A^{\otimes(m+2)}\rightarrow A^{\otimes(m+2)}$ by
\[
d_{(-1)}=d\text{ }
\]
and
\[
d_{(m)}=d\otimes1^{\otimes\left(  m+1\right)  }+1\otimes d_{(m-1)},
\]
which in more familiar form is
\begin{equation}
d_{(m)}=\sum_{i=0}^{m+1}1^{\otimes i}\otimes d\otimes1^{\otimes(m-i+1)}%
.\label{d-sub-m}%
\end{equation}
It is easy to check that $d_{(*)}$ is a differential with respect to internal
degree and
\begin{equation}
\partial\circ d_{(*)}-d_{(*-1)}\circ\partial=0.\label{commute1}%
\end{equation}
We refer to the double complex $\{TA,d_{(*)},\partial\}$ as the
\textit{two-sided bar resolution of }$A;$ the differentials $d_{(*)}$ and
$\partial$ have respective bidegree $(1,0)$ and $(0,-1)$ (see Figure 1). This
resolution is acyclic with respect to $\partial\,$via the contracting homotopy
$s$ given above.

An isomorphic (and more familiar) construction appears in \cite{Cartan} but
with the Eilenberg-Mac Lane shift in dimension. This dimension shift
introduces a set of signs that give rise to a standard double complex whose
subdiagrams anticommute; in this case $D=d_{(\ast)}+\partial$ is a
differential$.$ But short of that, we are better off without the dimension
shift since the subsequent signs unnecessarily complicate the exposition and
formulas. When total differentials are required, it is a simple matter to
introduce artificial signs; this is the strategy we adopt.\bigskip%
\[%
\begin{array}
[c]{c}%
\begin{array}
[c]{cccccccc}
& \;\vdots &  & \vdots &  & \vdots &  & \\
& {\scriptstyle\partial}\downarrow\; &  & \downarrow &  & \downarrow &  & \\
& \mathbf{k}^{\otimes4} & \rightarrow & (A^{\otimes4})^{1} & \rightarrow &
(A^{\otimes4})^{2} & \rightarrow & \cdots\\
& {\scriptstyle\partial}\downarrow\; &  & \downarrow &  & \downarrow &  & \\
& \mathbf{k}^{\otimes3} & \rightarrow & (A^{\otimes3})^{1} & \rightarrow &
(A^{\otimes3})^{2} & \rightarrow & \cdots\\
& {\scriptstyle\partial}\downarrow\; &  & \downarrow &  & \downarrow &  & \\
& \mathbf{k}^{\otimes2} & \rightarrow & (A^{\otimes2})^{1} & \rightarrow &
(A^{\otimes2})^{2} & \rightarrow & \cdots\\
& {\scriptstyle\partial}\downarrow\; &  & \downarrow &  & \downarrow &  & \\
& \mathbf{k} & \underset{d}{\rightarrow} & A^{1} & \underset{d}{\rightarrow} &
A^{2} & \underset{d}{\rightarrow} & \ldots
\end{array}
\\
\\
\text{The Bar Resolution}\\
\text{Figure 1.}%
\end{array}
\]

\subsection{The Cobar Resolution}

Now consider a \textit{d.g.c}. $(C,d,\Delta).$ The 2-sided cobar resolution of
$C$ is an $C$-free resolution dual to the 2-sided bar resolution.

As in (\ref{d-sub-m}), consider the differentials with respect to internal
degree $d_{(n)}:C^{\otimes\left(  n+2\right)  }\rightarrow C^{\otimes(n+2)}$
given by
\begin{equation}
\label{d-sub-n}d_{(n)}=\sum_{i=0}^{n+1}1^{\otimes i}\otimes d\otimes
1^{\otimes(n-i+1)}.
\end{equation}
Inductively define maps $\delta_{(n)}:C^{\otimes(n+2)}\rightarrow
C^{\otimes(n+3)}$ by
\[
\delta_{(-1)}=\Delta
\]
and
\[
\delta_{(n)}=\Delta\otimes1^{\otimes(n+1)}-1\otimes\delta_{(n-1)},
\]
which in more familiar form is
\begin{equation}
\label{cobar}\delta_{(n)}=\sum_{i=0}^{n+1}(-1)^{i}(1^{\otimes i}\otimes
\Delta\otimes1^{\otimes(n-i+1)}).
\end{equation}
Let $\delta=\sum_{n\geq-1}\delta_{(n)};$ using induction and the fact that
$\Delta$ is coassociative, it is easy to check that $\delta$ is a differential
with respect to external degree and%

\begin{equation}
\delta\circ d_{(\ast)}-d_{(\ast+1)}\circ\delta=0.\label{commute2}%
\end{equation}
We refer to the double complex $\{TC,d_{(\ast)},\delta\}$ as the
\textit{two-sided cobar resolution of }$C$. This resolution is acyclic with
respect to $\delta$ via contracting homotopy $\tau=\sum_{n\geq0}\tau_{n}$
where $\tau_{n}=[j_{1}\circ(\varepsilon\otimes1)]\otimes1^{\otimes n}.$

\section{Two-sided Differential Graded $\mathbf{k}$-Modules\ }

\subsection{Differential Graded $A$-Bimodules}

Let $\mathbf{k}$ be a field.

\begin{definition}
Let $(A,\mu)$ be a $\mathbf{k}$-algebra and let $M$ be a $\mathbf{k}$-module
for which there exist $\mathbf{k}$-linear structure maps $\lambda:A\otimes
M\rightarrow M$ and $\rho:M\otimes A\rightarrow M$ of degree zero such that

\begin{enumerate}
\item $\lambda\circ(\mu\otimes1)=\lambda\circ(1\otimes\lambda),$

\item $\rho\circ(1\otimes\mu)=\rho\circ(\rho\otimes1),$ and

\item $\lambda\circ(\eta\otimes1)\circ i_{1}=\rho\circ(1\otimes\eta)\circ
i_{2}=1$.
\end{enumerate}

\noindent Then $(M,\lambda,\rho)$ is a $A$\text{-bimodule\textit{; i}}t is
symmetric if $\rho=\lambda\circ(1,2).$ If $(M,\lambda,\rho)$ and $(M^{\prime
},\lambda^{\prime},\rho^{\prime})$ are $A$-bimodules, a map $f\in
Hom_{\mathbf{k}}^{*}(M,M^{\prime})$ is a map of $A$-bimodules if
$f\circ\lambda=\lambda^{\prime}\circ(1\otimes f)$ and $f\circ\rho=\rho
^{\prime}\circ(f\otimes1)$. The category of $A$-bimodules and $A$-bimodule
maps is denoted by $A$\textit{-bimod}.
\end{definition}

\begin{example}
Let $V$ be any $\mathbf{k}$-module, let $M=A\otimes V\otimes A,$ and consider
structure maps $\lambda^{\mu}=\mu\otimes1\otimes1$ and $\rho^{\mu}%
=1\otimes1\otimes\mu.$ Then $(A\otimes V\otimes A,\lambda^{\mu},\rho^{\mu})$
is an exterior $A$-bimodule; $\lambda^{\mu}$ and $\rho^{\mu}$ are called
exterior bimodule structure maps. This is not to be confused with the notion
of an exterior algebra.
\end{example}

Let $V$ be a $\mathbf{k}$-module, let $(M^{\prime},\lambda^{\prime}%
,\rho^{\prime})$ be an $A$-bimodule, and consider the exterior $A$-bimodule
$(A\otimes V\otimes A,\lambda^{\mu},\rho^{\mu}).$ There is a $\mathbf{k}%
$-linear isomorphism
\begin{equation}
\label{phi}\Phi_{V}:Hom_{A-bimod}^{*}(A\otimes V\otimes A,M^{\prime})\approx
Hom_{\mathbf{k}}^{*}(V,M^{\prime})
\end{equation}
given by $\Phi_{V}(f)=f\circ(\eta\otimes1\otimes\eta)\circ(i_{1}\otimes1)\circ
i_{2}$.

Let $W$ be a $\mathbf{k}$-module and consider the exterior $A$-bimodules
$A\otimes V\otimes A$ and $A\otimes W\otimes A.$ An $A$-bimodule map
$\Theta:A\otimes V\otimes A\rightarrow A\otimes W\otimes A$ of degree $p$
induces a $\mathbf{k}$-linear map
\begin{equation}
\label{THETA-star}\Theta^{*}:Hom_{\mathbf{k}}^{*}(W,M^{\prime})\rightarrow
Hom_{\mathbf{k}}^{*+p}(V,M^{\prime})
\end{equation}
given by $\Theta^{*}=\Phi_{V}\circ Hom_{A-bimod}(\Theta,M^{\prime})\circ
\Phi_{W}^{-1},$ where $\Phi_{W}^{-1}(g)=\lambda^{\prime}\circ(1\otimes
\rho^{\prime})\circ(1\otimes g\otimes1).$ A critical point here is the fact
that $\lambda^{\prime}\circ(1\otimes\rho^{\prime})\circ(1\otimes
g\otimes1):A\otimes V\otimes A\rightarrow M^{\prime}$ is an $A$-bimodule map;
the reader may wish to supply the proof.

\begin{definition}
Let $(A,d,\mu)$ be a d.g.a. and let $(M,\lambda,\rho)$ be an $A$-bimodule
equipped with a differential $d_{M}.$ Then $(M,d_{M})$ is a differential
graded (d.g.) $A$-bimodule\newline provided that

\begin{enumerate}
\item $d_{M}\circ\lambda=\lambda\circ(d\otimes1+1\otimes d_{M})$ and

\item $d_{M}\circ\rho=\rho\circ(d_{M}\otimes1+1\otimes d).$
\end{enumerate}

\noindent If $(M,d_{M})$ and $(M^{\prime},d_{M^{\prime}})$ are d.g.
$A$-bimodules, a map $f\in Hom_{A-bimod}^{*}(M,M^{\prime})$ is a map of
d.g.\textit{\ }$A$-bimodules if $(-1)^{\left|  f\right|  }f\circ
d_{M}-d_{M^{\prime}}\circ f=0$.
\end{definition}

\begin{example}
Every d.g.a. $(A,d,\mu)$ is a d.g. $A$-bimodule with respect to structure maps
$\lambda=\rho=\mu.$
\end{example}

\begin{example}
\label{external-d.g.bimod}Let $(A,d,\mu)$ be a d.g.a., let $m\geq0$, and
identify $A\otimes\mathbf{k}\otimes A$ with $A\otimes A.$ Consider the
differentials $d_{(m)}:A^{\otimes(m+2)}\rightarrow A^{\otimes(m+2)}$ and
$\partial_{(m)}:A^{\otimes(m+2)}\rightarrow A^{\otimes(m+1)}$ defined in
(\ref{d-sub-m}) and (\ref{bar}), respectively. Then $(A^{\otimes(m+2)}%
,d_{(m)})$ is an exterior d.g. $A$-bimodule, and by (\ref{commute1})$,$
$\partial_{(m)}$ is a map of exterior d.g.\textit{\ }$A$-bimodules.
\end{example}

Given \textit{d.g. }$A$-bimodules $(M,\lambda,\rho,d_{M})$ and $(M^{\prime
},\lambda^{\prime},\rho^{\prime},d_{M^{\prime}}),$ define a map\newline%
$\overline{d}:Hom_{A-bimod}^{*}(M,M^{\prime})\rightarrow Hom_{\mathbf{k}%
}^{*+1}(M,M^{\prime})$ by
\begin{equation}
\label{d-overbar}\overline{d}(f)=(-1)^{\left|  f\right|  }f\circ
d_{M}-d_{M^{\prime}}\circ f.
\end{equation}
The following fact will be useful in the construction that follows:

\begin{proposition}
\label{bimod-map}$\overline{d}(f)$ is a map of \textit{d.g. }$A$-bimodules.
\end{proposition}

\textit{Proof: }We check the compatibility of $\overline{d}(f)$ with the
structure map $\rho;$ the compatibility with $\lambda$ is similar.
$\overline{d}(f)\circ\rho=(-1)^{|f|}f\circ d_{M}\circ\rho-d_{M^{\prime}}\circ
f\circ\rho=(-1)^{|f|}f\circ\rho\circ(d_{M}\otimes1+1\otimes d)-d_{M^{\prime}%
}\circ f\circ\rho=(-1)^{|f|}f\circ\rho\circ(d_{M}\otimes1)+(-1)^{|f|}%
\rho^{\prime}\circ(f\otimes1)\circ(1\otimes d)-d_{M^{\prime}}\circ\rho
^{\prime}\circ(f\otimes1)=(-1)^{|f|}\rho^{\prime}\circ(f\otimes1)\circ
(d_{M}\otimes1+1\otimes d)-\rho^{\prime}\circ(d_{M^{\prime}}\otimes1+1\otimes
d)\circ(f\otimes1)=(-1)^{|f|}\rho^{\prime}\circ(f\circ d_{M}\otimes1)-$
$\rho^{\prime}\circ(d_{M^{\prime}}\circ f\otimes1)=\rho^{\prime}%
\circ(\overline{d}(f)\otimes1).$ It is trivial to check that $\overline
{d}(\overline{d}(f))=0; $ hence $\overline{d}(f)$ respects
differentials.\bigskip\ 

Let $(H,\mu,\Delta)$ be a $\mathbf{k}$-Hopf algebra.

\begin{definition}
Let $(M,\lambda,\rho)$ and $(M^{\prime},\lambda^{\prime},\rho^{\prime})$ be
$H$-bimodules. The internal (bimodule) tensor product of $(M,\lambda,\rho)$
with $(M^{\prime},\lambda^{\prime},\rho^{\prime})$ is the so called interior
$H$-bimodule $M\overline{\otimes}M^{\prime}=(M\otimes M^{\prime}%
,\lambda\overline{\otimes}\lambda^{\prime},\rho\overline{\otimes}\rho^{\prime
})$ with
\[
\lambda\overline{\otimes}\lambda^{\prime}=(\lambda\otimes\lambda^{\prime
})\circ(2,3)\circ(\Delta\otimes1\otimes1)
\]
and
\[
\rho\overline{\otimes}\rho^{\prime}=(\rho\otimes\rho^{\prime})\circ
(2,3)\circ(1\otimes1\otimes\Delta).
\]

\end{definition}

Since $\Delta$ is coassociative, $(\lambda\overline{\otimes}\lambda^{\prime
})\overline{\otimes}\lambda^{\prime\prime}=$ $\lambda\overline{\otimes
}(\lambda^{\prime}\overline{\otimes}\lambda^{\prime\prime})$ and
$(\rho\overline{\otimes}\rho^{\prime})\overline{\otimes}\rho^{\prime\prime}=$
$\rho$$\overline{\otimes}(\rho^{\prime}\overline{\otimes}\rho^{\prime\prime
}).$ Thus, the internal tensor product can be associatively applied to any
finite family of $H$-bimodules\textit{.}

\begin{definition}
Let $(M,\lambda,\rho)$ be an $H$-bimodule. The structure maps $\overline
{\lambda}^{\mu}$$=\mu\overline{\otimes}\lambda\overline{\otimes}\mu$ and
$\overline{\rho}^{\mu}$$=\mu\overline{\otimes}\rho\overline{\otimes}\mu$ on
the interior $H$-bimodule $H\overline{\otimes}M\overline{\otimes}H$ are called
the two-sided interior extensions of $\lambda$ and $\rho$ by $\mu$,
respectively$.$
\end{definition}

\begin{definition}
Let $(M,\lambda,\rho)$ be an $H$-bimodule. The interior $H$-bimodule
$M^{\overline{\otimes}n}=(M^{\otimes n},\lambda^{n},\rho^{n}),$ with
$\lambda^{n}=\lambda\overline{\otimes}\lambda^{n-1}=(\lambda\otimes
\lambda^{n-1})\circ(2,3)\circ(\Delta\otimes1^{\otimes n})$ and $\rho^{n}%
=\rho^{n-1}\overline{\otimes}\rho=\newline(\rho^{n-1}\otimes\rho
)\circ(n,n+1)\circ(1^{\otimes n}\otimes\Delta),$ is called the n-fold interior
(bimodule) tensor power of $M$.
\end{definition}

\begin{example}
For each $n\geq1,$ the n-fold interior (bimodule) tensor power of $H$ is the
interior $H$-bimodule $H^{\overline{\otimes}n}=(H^{\otimes n},\lambda^{n}%
,\rho^{n})$ with
\[
\lambda^{n}=\mu^{\otimes n}\circ(1\ 3\ 5\ \cdots\ (2n-1)\ 2\ 4\ 6\ \cdots
\ 2n)\circ\ \prod_{i=n}^{2n-2}(\Delta\otimes1^{\otimes(3n-i-2)})\ \text{and}
\]
\[
\rho^{n}=\mu^{\otimes n}\circ(1\ 3\ 5\ \cdots\ (2n-1)\ 2\ 4\ 6\ \cdots
\ 2n)\circ\ \prod_{i=n}^{2n-2}(1^{\otimes(3n-i-2)}\otimes\Delta).\
\]

\end{example}

\begin{example}
Let $(H,d,\mu,\Delta)$ be a d.g.h.a., let $n\geq0,$ and identify the interior
$H$-bimodules $H\overline{\otimes}\mathbf{k}\overline{\otimes}H$ and
$H\overline{\otimes}H.$ Then $(H^{\overline{\otimes}(n+2)},d_{(n)})$ is an
interior d.g. $H$-bimodule and the differential $\delta_{(n)}:H^{\overline
{\otimes}(n+2)}\rightarrow H^{\overline{\otimes}(n+3)}$ defined in
(\ref{cobar}) is a d.g. $H$-bimodule map.\bigskip
\end{example}

\subsection{Differential Graded $C$-Bicomodules}

Let $\mathbf{k}$ be a field.

\begin{definition}
Let $(C,\Delta)$ be a $\mathbf{k}$-coalgebra and let $N$ be a $\mathbf{k}%
$-module for which there exist $\mathbf{k}$-linear structure maps
$\lambda:N\rightarrow C\otimes N$ and $\rho:N\rightarrow N\otimes C$ of degree
zero such that

\begin{enumerate}
\item $(\Delta\otimes1)\circ\lambda=(1\otimes\lambda)\circ\lambda,$

\item $(1\otimes\Delta)\circ\rho=(\rho\otimes1)\circ\rho,$ and

\item $j_{1}\circ(\varepsilon\otimes1)\circ\lambda=j_{2}\circ(1\otimes
\varepsilon)\circ\rho=1$.
\end{enumerate}

\noindent Then the triple $(N,\lambda,\rho)$ is an $C$%
\text{\textit{-bicomodule.}} If $(N,\lambda,\rho)$ and $(N^{\prime}%
,\lambda^{\prime},\rho^{\prime})$ are $C$-\newline bicomodules, a map $g\in
Hom_{\mathbf{k}}^{*}(N,N^{\prime})$ is a map of $C$-bicomodules if
$\lambda^{\prime}\circ g=(1\otimes g)\circ\lambda$ and $\rho^{\prime}\circ
g=(g\otimes1)\circ\rho$. The category of $C$-bicomodules and $C$-bicomodule
maps is denoted by $C$-\textit{bicomod.}
\end{definition}

\begin{example}
Let $V$ be any $\mathbf{k}$-module, let $N=C\otimes V\otimes C,$ and consider
structure maps $\lambda_{\Delta}=\Delta\otimes1\otimes1$ and $\rho_{\Delta
}=1\otimes1\otimes\Delta.$ Then $(C\otimes V\otimes C,\lambda_{\Delta}%
,\rho_{\Delta})$ is an exterior $C$-bicomodule; $\lambda_{\Delta}$ and
$\rho_{\Delta}$ are called exterior bicomodule structure maps.
\end{example}

Let $V$ be any $\mathbf{k}$-module, let $(N,\lambda,\rho)$ be any
$C$-bicomodule, and consider the exterior $C$-bicomodule $(C\otimes V\otimes
C,\lambda_{\Delta},\rho_{\Delta}).$ There is a $\mathbf{k}$-linear
isomorphism
\begin{equation}
\label{psi}\Psi_{V}:Hom_{C-bicomod}^{*}(N,C\otimes V\otimes C)\approx
Hom_{\mathbf{k}}^{*}(N,V)
\end{equation}
given by $\Psi_{V}(g)=j_{1}\circ(1\otimes j_{2})\circ(\varepsilon
\otimes1\otimes\varepsilon)\circ g$.

Let $W$ be a $\mathbf{k}$-module and consider the exterior $C$-bicomodules
$C\otimes V\otimes C$ and $C\otimes W\otimes C.$ An $C$-bicomodule map $\Xi
$$:C\otimes V\otimes C\rightarrow C\otimes W\otimes C$ of degree $q$ induces a
$\mathbf{k}$-linear map
\begin{equation}
\label{XI-star}\Xi_{*}:Hom_{\mathbf{k}}^{*}(N,V)\rightarrow Hom_{\mathbf{k}%
}^{*+q}(N,W)
\end{equation}
given by $\Xi_{*}=\Psi_{W}\circ Hom_{C-bicomod}(N,\Xi)\circ\Psi_{V}^{-1},$
where $\Psi_{V}^{-1}(g)=(1\otimes g\otimes1)\circ(1\otimes\rho)\circ\lambda.$
Again, it is important to note that $(1\otimes g\otimes1)\circ(1\otimes
\rho)\circ\lambda:N\rightarrow C\otimes V\otimes C$ is an $C$-bicomodule map.\ 

\begin{definition}
Let $(C,d,\Delta)$ be a d.g.c. and let $(N,\lambda,\rho)$ be an $C$%
-bicomodule\newline equipped with a differential $d_{N}.$ Then $(N,d_{N})$ is
a d.g. $C $-bicomodule provided that
\end{definition}

\begin{enumerate}
\item $\lambda\circ d_{N}=(d\otimes1+1\otimes d_{N})\circ\lambda$ and

\item $\rho\circ d_{N}=(d_{N}\otimes1+1\otimes d)\circ\rho.$
\end{enumerate}

\noindent If $(N,d_{N})$ and $(N^{\prime},d_{N^{\prime}})$ are d.g.
$C$-bicomodules, a map $g\in Hom_{C-bicomod}^{*}(N,N^{\prime})$ is a map of
d.g. $C$-bicomodules if $(-1)^{|g|}g\circ d_{N}-d_{N^{\prime}}\circ g=0.$

\begin{example}
Every d.g.c. $(C,d,\Delta)$ is a d.g. $C$-bicomodule with respect to structure
maps $\lambda=\rho=\Delta.$
\end{example}

\begin{example}
Let $(C,d,\Delta)$ be a d.g.c., let $n\geq0$, and identify $C\otimes
\mathbf{k}\otimes C$ with $C\otimes C.$ Consider the differentials
$d_{(n)}:C^{\otimes(n+2)}\rightarrow C^{\otimes(n+2)}$ and $\delta
_{(n)}:C^{\otimes(n+2)}\rightarrow C^{\otimes(n+3)}$ defined in (\ref{d-sub-n}%
) and (\ref{cobar}), respectively. Then $(C^{\otimes(n+2)},d_{(n)}) $ is an
exterior d.g. $C$-bicomodule, and by (\ref{commute2}), $\delta_{(n)} $ is a
map of exterior \textit{d.g. }$C$-bimodules.
\end{example}

Given \textit{d.g. }$C$-bicomodules $(N,\lambda,\rho,d_{N})$ and $(N^{\prime
},\lambda^{\prime},\rho^{\prime},d_{N^{\prime}}),$ define a map $\overline
{d}:Hom_{C-bicomod}^{*}(N,N^{\prime})\rightarrow Hom_{\mathbf{k}}%
^{*+1}(N,N^{\prime})$ as in (\ref{d-overbar}). The reader can check that:

\begin{proposition}
\label{bicomod-map}$\overline{d}(g)$ is a map of d.g.\textit{\ }$C$-bicomodules.
\end{proposition}

Let $(H,\mu,\Delta)$ be a $\mathbf{k}$-Hopf algebra$.$\ 

\begin{definition}
Let $(N,\lambda,\rho)$ and $(N^{\prime},\lambda^{\prime},\rho^{\prime})$ be
$H$-bicomodules. The internal (bicomodule) tensor product of $(N,\lambda
,\rho)$ with $(N^{\prime},\lambda^{\prime},\rho^{\prime})$ is the so called
interior $H$-bicomodule $N\underline{\otimes}N^{\prime}=(N\otimes N^{\prime
},\lambda\underline{\otimes}\lambda^{\prime},\rho\underline{\otimes}%
\rho^{\prime})$ with
\[
\lambda\underline{\otimes}\lambda^{\prime}=(\mu\otimes1\otimes1)\circ
(2,3)\circ(\lambda\otimes\lambda^{\prime})
\]
and
\[
\rho\underline{\otimes}\rho^{\prime}=(1\otimes1\otimes\mu)\circ(2,3)\circ
(\rho\otimes\rho^{\prime}).
\]

\end{definition}

Since $\mu$ is associative, $(\lambda\underline{\otimes}\lambda^{\prime
})\underline{\otimes}\lambda^{\prime\prime}=\lambda\underline{\otimes}%
(\lambda^{\prime}\underline{\otimes}\lambda^{\prime\prime})$ and
$(\rho\underline{\otimes}\rho^{\prime})\underline{\otimes}\rho^{\prime\prime
}=$ $\rho\underline{\otimes}(\rho^{\prime}\underline{\otimes}\rho
^{\prime\prime})$. Thus, the internal tensor product can be associatively
applied to any finite family of $H$-bicomodules\textit{.}

\begin{definition}
Let $(N,\lambda,\rho)$ be an $H$-bicomodule. The structure maps $\overline
{\lambda}_{\Delta}=\Delta\underline{\otimes}\lambda\underline{\otimes}\Delta$
and $\overline{\rho}_{\Delta}=\Delta\underline{\otimes}\rho\underline{\otimes
}\Delta$ on the interior $H$-bicomodule $H\underline{\otimes}N\underline
{\otimes}H$ are called the two-sided interior extensions of $\lambda$ and
$\rho$ by $\Delta,$ respectively.
\end{definition}

\begin{definition}
Let $(N,\lambda,\rho)$ be an $H$-bicomodule. The interior $H$-bicomodule
$N^{\underline{\otimes}m}=(N^{\otimes m},\lambda_{m,}\rho_{m}),$ with
$\lambda_{m}=$ $\lambda\underline{\otimes}\lambda_{m-1}=(\mu\otimes1^{\otimes
m})\circ(2,3)\circ(\lambda\otimes\lambda_{m-1})$ and $\rho_{m}=\rho
_{m-1}\underline{\otimes}\rho=\newline(1^{\otimes m}\otimes\mu)\circ
(m,m+1)\circ(\rho_{m-1}\otimes\rho),$ is called the $m$-fold internal
(bicomodule) tensor power of $N.$
\end{definition}

\begin{example}
Consider the $H$-bicomodule $(H,\lambda,\rho)$ where $\lambda=\rho=\Delta.$
For each $m\geq1,$ the $m$-fold internal (bicomodule) tensor power of $H$ is
the interior $H$-bicomodule $H^{\underline{\otimes}m}=(H^{\otimes m}%
,\lambda_{m},\rho_{m})$ with
\[
\lambda_{m}=\prod_{i=m}^{2m-2}(\mu\otimes1^{\otimes i})\circ(1\ 3\ 5\ \cdots
\ (2m-1)\ 2\ 4\ 6\ \cdots\ 2m)^{-1}\circ\Delta^{\otimes m}\ \text{ and}
\]
\[
\rho_{m}=\prod_{i=m}^{2m-2}(1^{\otimes i}\otimes\mu)\circ(1\ 3\ 5\ \cdots
\ (2m-1)\ 2\ 4\ 6\ \cdots\ 2m)^{-1}\circ\Delta^{\otimes m}\ .
\]

\end{example}

\begin{example}
Let $(H,d,\mu,\Delta)$ be a d.g.h.a., let $m\geq0,$ and identify the interior
$H$-bicomodules $H\underline{\otimes}\mathbf{k}\underline{\otimes}H$ and
$H\underline{\otimes}H.$ Then $(H^{\underline{\otimes}(m+2)},d_{(m)})$ is an
interior d.g. $H$-bicomodule$.$ Furthermore, $\partial_{(m)}:H^{\underline
{\otimes}(m+2)}\rightarrow$ $H^{\underline{\otimes}(m+1)}$ defined in
(\ref{bar}) is a map of interior d.g. $H$-bicomodules.\bigskip\ 
\end{example}

\subsection{Differential Graded $H$-bidimodules}

Let $\mathbf{k}$ be a field and let $(H,\mu,\Delta)$ be a $\mathbf{k}$-Hopf algebra.

\begin{definition}
\label{bidimod}Let $E$ be a $\mathbf{k}$-module such that $(E,\lambda
^{\#},\rho^{\#})$ is an $H$-bimodule and $(E,\lambda_{\#},\rho_{\#})$ is an
$H$-bicomodule. Then $(E,\lambda^{\#},\rho^{\#},\lambda_{\#},\rho_{\#})$ is
an\textit{\ }$H$-bidimodule if its $H$-bimodule and $H$-bicomodule structures
are compatible in the following sense:

\begin{enumerate}
\item $\lambda^{\#}\in Hom_{H-bicomod}(H\underline{\otimes}E,E),$

\item $\rho^{\#}\in Hom_{H-bicomod}(E\underline{\otimes}H,E),$

\item $\lambda_{\#}\in Hom_{H\text{-}bimod}(E,H\overline{\otimes}E),$ and

\item $\rho_{\#}\in Hom_{H\text{-}bimod}($$E,E\overline{\otimes}H)$.
\end{enumerate}

\noindent A map of $H$-bidimodules preserves bimodule and bicomodule
structure. Denote the category of $H$-bidimodules and $H$-bidimodule maps by
\textit{$H$-bidimod.}
\end{definition}

\begin{example}
\label{bi-dimod}Let $(M,\lambda,\rho)$ be any $H$-bimodule. Then
$H\overline{\otimes}M\overline{\otimes}H=(H\otimes M\otimes H,\overline
{\lambda}^{\mu},\overline{\rho}^{\mu},\lambda_{\Delta},\rho_{\Delta})$ is an
$H$-bidimodule, where $(H\otimes M\otimes H,\overline{\lambda}^{\mu}%
,\overline{\rho}^{\mu})$ is an interior $H$-bimodule and $(H\otimes M\otimes
H,\lambda_{\Delta},\rho_{\Delta})$ is an exterior $H$-bicomodule. Dually, if
$(N,\lambda,\rho)$ is any $H$-bicomodule, then $H\underline{\otimes
}N\underline{\otimes}H=(H\otimes N\otimes H,\lambda^{\mu},\rho^{\mu}%
,\overline{\lambda}_{\Delta},\overline{\rho}_{\Delta})$ is an $H$-bidimodule.
In particular:
\end{example}

\begin{example}
$H^{\underline{\otimes}(m+2)}$ and $H^{\overline{\otimes}(n+2)}$ are $H$-bidimodules.
\end{example}

Isomorphisms (\ref{phi}) and (\ref{psi}) extend to the $H$-bidimodules in
Example \ref{bi-dimod} giving
\[
\Phi_{N}:Hom_{H-bidimod}^{*}(H\underline{\otimes}N\underline{\otimes
}H,H\overline{\otimes}M\overline{\otimes}H)\approx Hom_{H-bicomod}%
^{*}(N,H\overline{\otimes}M\overline{\otimes}H)
\]
and
\[
\Psi_{M}:Hom_{H-bidimod}^{*}(H\underline{\otimes}N\underline{\otimes
}H,H\overline{\otimes}M\overline{\otimes}H)\approx Hom_{H-bimod}%
^{*}(H\underline{\otimes}N\underline{\otimes}H,M).
\]
Thus,
\[
\Phi_{N}\circ\Psi_{M}=\Psi_{M}\circ\Phi_{N}:Hom_{H-bidimod}^{*}(H\underline
{\otimes}N\underline{\otimes}H,H\overline{\otimes}M\overline{\otimes}H)\approx
Hom_{\mathbf{k}}^{*}(N,M).
\]
In particular, for $m,n\geq0,$ $N=H^{\underline{\otimes}m}$ and
$M=H^{\overline{\otimes}n}$ we have
\begin{equation}
\Phi_{N}\circ\Psi_{M}:Hom_{H-bidimod}^{*}(H^{\underline{\otimes}%
(m+2)},H^{\overline{\otimes}(n+2)})\approx Hom_{\mathbf{k}}^{*}(H^{\underline
{\otimes}m},H^{\overline{\otimes}n}).\label{isobidimod}%
\end{equation}

Now consider the $H$-bidimodules $H\underline{\otimes}N\underline{\otimes}H,$
$H\underline{\otimes}N^{\prime}\underline{\otimes}H$ and $H\overline{\otimes
}M\overline{\otimes}H.$ If $\Theta:\newline H\underline{\otimes}%
N\underline{\otimes}H\rightarrow H\underline{\otimes}N^{\prime}\underline
{\otimes}H$ is an $H$-bidimodule map of degree $p,$ there is an induced map
$\Theta^{*}:Hom_{\mathbf{k}}^{*}(N^{\prime},M)\rightarrow Hom_{\mathbf{k}%
}^{*+p}(N,M)$ given by
\begin{equation}
\label{THETA-STAR}\Theta^{*}=(\Phi_{N}\circ\Psi_{M})\circ Hom_{H-bidimod}%
^{*}(\Theta,H\overline{\otimes}M\overline{\otimes}H)\circ(\Phi_{N^{\prime}%
}\circ\Psi_{M})^{-1}.
\end{equation}
Dually, given $H$-bidimodules $H\overline{\otimes}M\overline{\otimes}H,$
$H\overline{\otimes}M^{\prime}\overline{\otimes}H,$ and $H\underline{\otimes
}N\underline{\otimes}H,$ an $H$-bidimodule map $\Xi:H\overline{\otimes
}M\overline{\otimes}H\rightarrow H\overline{\otimes}M^{\prime}\overline
{\otimes}H$ of degree $q$ induces a map $\Xi_{*}:Hom_{\mathbf{k}}%
^{*}(N,M)\rightarrow Hom_{\mathbf{k}}^{*+q}(N,M^{\prime})$ via
\begin{equation}
\label{XI-STAR}\Xi_{*}=(\Psi_{M^{\prime}}\circ\Phi_{N})\circ Hom_{H-bidimod}%
^{*}(H\underline{\otimes}N\underline{\otimes}H,\Xi)\circ(\Psi_{M}\circ\Phi
_{N})^{-1}.
\end{equation}

\begin{definition}
Let $(E,\lambda^{\#},\rho^{\#},\lambda_{\#},\rho_{\#})$ be an $H$-bidimodule
equipped with a differential $d_{E}.$ Then $(E,d_{E})$ is a d.g.
$H$-bidimodule if $(E,\lambda^{\#},\rho^{\#},d_{E})$ is a d.g. $H$-bimodule
and $(E,\lambda_{\#},\rho_{\#},d_{E})$ is a d.g. $H$-bicomodule. If
$(E,d_{E})$ and $(E^{\prime},d_{E}^{\prime})$ are d.g. $H$-bidimodules, a map
$h\in Hom_{H-bidimod}^{*}(E,E^{\prime})$ is a d.g. $H$-bidimodule map if
$(-1)^{\left|  h\right|  }h\circ d_{E}-d_{E^{\prime}}\circ h=0$.
\end{definition}

\begin{example}
\label{d.g.bidimod}For each $m,n\geq0,$ $(H^{\underline{\otimes}(m+2)}%
,d_{(m)})$ and $(H^{\overline{\otimes}(n+2)},d_{(n)})$ are d.g. $H$%
-bidimodules; maps $\partial_{(m)}$ and $\delta_{(n)}$ are d.g. $H$-bidimodule maps.
\end{example}

Let $(E,\lambda^{\#},\rho^{\#},\lambda_{\#},\rho_{\#},d_{E})$ and $(E^{\prime
},\lambda^{\#^{\prime}},\rho^{\#^{\prime}},\lambda_{\#}^{\prime},\rho
_{\#}^{\prime},d_{E}^{\prime})$ be d.g. $H$-bidimodules and define
$\overline{d}:Hom_{H-bidimod}^{*}(E,E^{\prime})\rightarrow Hom_{\mathbf{k}%
}^{*+1}(E,E^{\prime})$ by $\overline{d}(h)=(-1)^{\left|  h\right|  }h\circ
d_{E}-d_{E^{\prime}}\circ h$.

\begin{proposition}
\label{bidimod-map}Then $\overline{d}(h)$ is a map of d.g. $H$-bidimodules.
\end{proposition}

\section{The Deformation Complex for Differential Graded Structures}

\subsection{The Deformation Complex for Differential Graded Algebras}

We begin by defining the Hochschild cohomology with coefficients in a
\textit{d.g.} $A$-bimodule. Let $(A,d,\mu)$ be a \textit{d.g.a. }and let
$(M,\lambda,\rho,d_{M})$ be a \textit{d.g.} $A$-bimodule. For each $m\geq0,$
consider the exterior \textit{d.g. }$A$-bimodule $(A^{\otimes(m+2)},d_{(m)})$,
where $d_{(m)}$ is defined as in (\ref{d-sub-m}). By Proposition
\ref{bimod-map}, there is a map $d^{p,m}:Hom_{A-bimod}^{p}(A^{\otimes
(m+2)},M)\rightarrow Hom_{A-bimod}^{p+1}(A^{\otimes(m+2)},M)$ given by
\[
d^{p,m}(f)=(-1)^{p}f\circ d_{(m)}-d_{M}\circ f.
\]
It is easy to check that $d^{*,*}$ is a differential of bidegree $(1,0)$.
Consider the $\mathbf{k}$-linear isomorphism $\Phi:Hom_{A-bimod}%
^{*}(A^{\otimes(m+2)},M)\rightarrow Hom_{\mathbf{k}}^{*}(A^{\otimes m},M)$
defined in (\ref{phi}). Then $d^{*,*}$ induces a $\mathbf{k}$-linear
differential $d_{B}^{p,m}:Hom_{\mathbf{k}}^{p}(A^{\otimes m},M)\rightarrow
Hom_{\mathbf{k}}^{p+1}(A^{\otimes m},M)$ of bidegree $(1,0)$ via
\[
d_{B}^{*,*}=\Phi\circ d^{*,*}\circ\Phi^{-1}.
\]
Furthermore, the $A$-bimodule map $\partial_{(m)}:A^{\otimes(m+2)}\rightarrow
A^{\otimes(m+1)}$ induces a map $\partial^{*,m}=Hom_{A-bimod}^{*}%
(\partial_{(m)},M)$ and subsequently, as in (\ref{THETA-star}), a $\mathbf{k}%
$-linear map $\partial_{B}^{*,*}:Hom_{\mathbf{k}}^{*}(A^{\otimes*}%
,M)\newline\rightarrow Hom_{\mathbf{k}}^{*}(A^{\otimes(*+1)},M)$ of bidegree
$(0,1)$ defined by
\[
\partial_{B}^{*,*}=\Phi\circ\partial^{*,*+1}\circ\Phi^{-1}.
\]
It is easy to check that $\partial_{B}$ is a differential; the fact that
$\partial_{B}\circ d_{B}-d_{B}\circ\partial_{B}=0$ follows easily from
(\ref{commute1})$.$ At a particular $f\in Hom_{\mathbf{k}}^{p}(A^{\otimes
m},M),$ the differentials $d_{B}$ and $\partial_{B}$ can be written as
\begin{equation}
d_{B}^{p,m}(f)=(-1)^{p}f\circ d_{(m-2)}-d_{M}\circ f\text{ \ }\label{d}%
\end{equation}
and
\begin{equation}
\partial_{B}^{p,m}(f)=\lambda\circ(1\otimes f)-f\circ\partial_{(m-1)}%
+(-1)^{m+1}\rho\circ(f\otimes1),\label{d-bar}%
\end{equation}
where $d_{(-2)}=0$ and $\partial_{(-1)}=0.$

Refer to $B^{p,m}(A;M)=Hom_{\mathbf{k}}^{p}(A^{\otimes m},M)$ as the space of
\textit{Hochschild m-cochains on }$A$ \textit{of degree p}. The double complex
$\{B^{\ast,\ast}(A;M),d_{B},\partial_{B}\}$ is called the \textit{Hochschild
cochain complex on }$A$\textit{\ with coefficients in the d.g. }%
$A$\textit{-bimodule }$M$ (see Figure 2)\textit{. }Define the space of total
$r$-cochains by $B^{r}(A;M)=\newline\sum_{p\in\mathbb{Z};\text{ }r\geq
p}B^{p,r-p}(A;M)$ and define $D_{B}$ on the component $B^{p,r-p}(A;M)$ by
\[
D_{B}=d_{B}-(-1)^{p}\partial_{B}.
\]
Then $D_{B}$ is a total differential; the sign $-(-1)^{p}$ is introduced so
that $D_{B}^{2}=0.$ Now define the \textit{Hochschild cohomology of }$A$
\textit{with coefficients in }$M,$ denoted by $H_{d.g.a}^{\ast}(A;M),$ to be
the homology of the total complex $\{B^{\ast}(A;M),\ D_{B}\}.\,$\bigskip\ 

\noindent\hspace*{-0.3in}%
\[%
\begin{array}
[c]{c}%
\begin{array}
[c]{ccccccc}
& \;\;\ \vdots &  & \vdots &  & \vdots & \\
& {\scriptstyle\partial}_{\scriptscriptstyle B}\uparrow\; &  & \uparrow &  &
\uparrow & \\
&  &  &  &  &  & \\
\cdots & Hom_{\mathbf{k}}^{p-1}(A^{\otimes3},M) & \rightarrow &
Hom_{\mathbf{k}}^{p}(A^{\otimes3},M) & \rightarrow & Hom_{\mathbf{k}}%
^{p+1}(A^{\otimes3},M) & \cdots\\
&  &  &  &  &  & \\
& {\scriptstyle\partial}_{\scriptscriptstyle B}\uparrow\; &  & \uparrow &  &
\uparrow & \\
&  &  &  &  &  & \\
\cdots & Hom_{\mathbf{k}}^{p-1}(A^{\otimes2},M) & \rightarrow &
Hom_{\mathbf{k}}^{p}(A^{\otimes2},M) & \rightarrow & Hom_{\mathbf{k}}%
^{p+1}(A^{\otimes2},M) & \cdots\\
&  &  &  &  &  & \\
& {\scriptstyle\partial}_{\scriptscriptstyle B}\uparrow\; &  & \uparrow &  &
\uparrow & \\
&  &  &  &  &  & \\
\cdots & Hom_{\mathbf{k}}^{p-1}(A^{\otimes1},M) & \rightarrow &
Hom_{\mathbf{k}}^{p}(A^{\otimes1},M) & \rightarrow & Hom_{\mathbf{k}}%
^{p+1}(A^{\otimes1},M) & \cdots\\
&  &  &  &  &  & \\
\noalign {.\ .\ .\ .\ .\ .\ .\ .\ .\ .\ .\ .\ .\ .\ .\ .\ .\ .\ .\ .\ .\ .\ .\ .\ .\ .\ .\ .\ .\ .\ .\ .\ .\ .\ .\ .\ .\ .\ .\ .\ .\ .\ .\ .\ .\ .\ .\ .\ .\ .\ .\ .\ .\ .\ .} &
&  &  &  &  & \\
& {\scriptstyle\partial}_{\scriptscriptstyle B}\uparrow\; &  & \uparrow &  &
\uparrow & \\
&  &  &  &  &  & \\
\cdots & Hom_{\mathbf{k}}^{p-1}(\mathbf{k},M) & \rightarrow & Hom_{\mathbf{k}%
}^{p}(\mathbf{k},M) & \rightarrow & Hom_{\mathbf{k}}^{p+1}(\mathbf{k},M) &
\ldots\\
&  & d_{\scriptscriptstyle B} &  & d_{\scriptscriptstyle B} &  &
\end{array}
\\
\\
\text{The Hochschild cochain complex on\textit{\ }}A\\
\text{Figure 2.}%
\end{array}
\]

We say that an $m$-cochain $f\in B^{p,m}(A;M)$ is \textit{normalized} if
$f(a_{1}\otimes\cdots\otimes a_{m})=0$ whenever $\left|  a_{i}\right|  =0$ for
some $i=1,2,\ldots,m.$ The differential $\partial_{B}$ restricts to the
subspace of normalized cochains and a standard theorem \cite{MacLane} assures
that the subcomplex given by such a restriction is cochain homotopic to the
full cochain complex. Consequently, we shall use normalized cochains; the
symbol $B^{p,m}(A;M)$ will henceforth denote the space of normalized
$m$-cochains of degree $p$.

The applications require a bitruncation of the Hochschild cochain complex.
First, delete the bottom row in Figure 2 above; the subsequent theory is
referred to as the \textit{restricted }Hochschild cohomology of $A$. Denote
the space of restricted total cochains by $\widetilde{B}^{*}(A;M)$ and the
restricted cohomology by $\widetilde{H}^{*}(A;M).$ For $n\in\{3,4,\ldots
\}\cup\{\infty\},$ further restrict the Hochschild complex to those cochains
in bidegree $(p,m)$ with $p\geq3-n$ and $m\geq1.$ Denote the bitruncated total
$r$-cochains by $\widetilde{B}^{r}(A;M;n)=\sum_{\text{ }r-1\geq p\geq
3-n}B^{p,r-p}(A;M).$ The homology of the complex $\{\widetilde{B}%
^{*}(A;M;n),\ D_{B}\},$ denoted by $\widetilde{H}_{d.g.a.}^{*}(A;M;n)$, is
called the \textit{restricted Hochschild cohomology of }$A$ \textit{truncated
at degree }$3-n.$

The \textit{deformation complex for $A$ as a strict }$A\left(  n\right)
$\textit{-algebra }is the cochain complex $\{\widetilde{B}^{*}(A;A;n),D_{B}%
\}$\textit{. }The cohomology $\widetilde{H}_{d.g.a.}^{*}(A;A;n)$ directs the
deformation theory in the following way: If $A_{t}=(A[[t]],\ \mu_{t}%
^{(1)},\ \mu_{t}^{(2)},\ \ldots)$ is a deformation of $A$ as a strict
$A(n)$-algebra, we agree that for $1\leq i\leq n$ the maps $\mu_{t}^{(i)}%
=\mu_{0}^{(i)}+t\mu_{1}^{(i)}+t^{2}\mu_{2}^{(i)}+\cdots$ satisfy $\mu
_{1}^{(i)}\in B^{2-i,i}(A;A)\ $with $\mu_{0}^{(1)}=d,$ $\mu_{0}^{(2)}=\mu,$
and $\mu_{0}^{(i)}=0$ for $i>2.$ Then $a=\sum_{i=1}^{n}\mu_{1}^{(i)}$
$\in\widetilde{B}^{2}(A;A;n)$ satisfies $D_{B}(a)=0. $ Furthermore, given a
total cocycle $a=\sum_{i=1}^{n}\mu_{1}^{(i)}$ $\in\widetilde{B}^{2}(A;A;n),$
the obstructions to extending the corresponding linear approximation
$(d+t\mu_{1}^{(1)},$ $\mu+t\mu_{1}^{(2)},$ $t\mu_{1}^{\left(  3\right)  },$
$\ldots,\ t\mu_{1}^{(i)},\ldots)$ to a deformation appear as an inductively
defined sequence of cocycles in $\widetilde{B}^{3}(A;A;n).$ We note that the
deformation theory of $A$ as an $A(\infty)$-algebra also appears in
\cite{Kontsevich} and \cite{Penkava}. The case $n=3$ is discussed in some
detail in section 6 (as a special case) and subsequently by Lazarev and
Movshev in the sequel \cite{Lazarev}. When $n=3,$ we adopt the standard
notation $A_{t}=(A[[t]],\ d_{t},\ \mu_{t})$ with $d_{t}=d+td_{1}+t^{2}%
d_{2}+\cdots$ and $\mu_{t}=\mu+t\mu_{1}+t^{2}\mu_{2}+\cdots,$ in which case
$A_{t}$ is an (associative) \textit{d.g.a. } Now suppose that such an $A$ is
commutative\textit{.}

If $(A,d,\mu)$ is a commutative \textit{d.g.a. (c.d.g.a.),} a general
deformation $A_{t}=\newline(A[[t]],\ d_{t,}\ \mu_{t})$ fails to be
commutative. In some applications, such as the classification of rational
homotopy type for example, one desires only commutative deformations; in this
case it is necessary to restrict the Hochschild cochain complex to the
subcomplex of cochains with the potential to spawn commutative deformations. A
discussion of this subcomplex, called the \textit{Harrison cochain complex on
A, }now follows.

Let $(A,d,\mu)$ be a \textit{c.d.g.a.} and let $(M,\lambda,\rho,d_{M})$ be a
symmetric \textit{d.g. }$A$-bimodule. The desired cochains $f\in B^{*,2}(A;M)$
are the symmetric functions, and in general, the desired cochains $f\in
B^{*,m}(A;M)$ vanish on sums of certain shuffle permutations. Precisely, let
$\sigma_{r,s}\in S_{n}$ denote a $(r,s)$-shuffle \cite{MacLane}; denote its
sign as a permutation by $(-1)^{\sigma_{r,s}}.$ Then $f\in B^{p,m}(A;M)$ is a
\textit{Harrison }$m$\textit{-cochain on }$A$ \textit{in degree p} if and only
if $f $ vanishes on $\sum_{\sigma_{r,m-r}}(-1)^{\sigma_{r,m-r}}\sigma_{r,m-r}$
for each $r=1,2,\ldots,m-1$. The space of \textit{Harrison m-cochains of
degree p} is denoted by $Ch^{p,m}(A;M);$ note that $Ch^{*,1}(A;M)=B^{*,1}%
(A;M).$ The differentials $d_{B}$ and $\partial$$_{B}$ restrict $\,$to
$Ch^{*,*}(A;M)$, so let $Ch^{r}(A;M)=\sum_{p+m=r}Ch^{p,m}(A;M)$ and define the
\textit{Harrison cohomology of }$A,$ $Harr_{c.d.g.a.}^{*}(A;M),$ to be the
homology of the total complex $\{Ch^{*}(A;M),D_{B}\}.$ As in the general
Hochschild case, the symbol $\widetilde{C}h^{*}(A;M;n)$ denotes the total
bitruncated Harrison cochains and $\widetilde{H}arr_{c.d.g.a.}^{*}(A;M;n)$
denotes the corresponding cohomology. The complex $\{\widetilde{C}%
h^{*}(A;A;n),D_{B}\}$ is the deformation complex for $A$ as a ''balanced''
$A(n)$-algebra; see \cite{Kadeishvili} and \cite{Markl}.

If $(A,\mu)$ is a \textit{c.g.a.} sans differential, one can forget the
internal grading and grade the Harrison cohomology externally (with respect to
the number of tensor factors) as one does classically. Let $Harr^{*}(A;M)$
denote the Harrison cohomology graded in this way; one has the following
result, which is a consequence of the Hochschild, Kostant, and Rosenberg
Theorem \cite{Hochschild}:\ 

\begin{theorem}
\label{free}Let $\mathbf{k}$ be a field of characteristic $0$. If $A$ is a
free commutative $\mathbf{k}$-algebra and $(M,\lambda,\mu)$ is any symmetric
$A$-bimodule, then $Harr^{n}(A;M)=0$ whenever $n>1.$
\end{theorem}

\noindent The requirement that $\mathbf{k}$ have characteristic $0$ is
critical here; the result fails in characteristic $p>0$ \cite{Barr}$.$ Theorem
\ref{free} allows us to view the rational perturbation theory of Felix and
Halperin-Stasheff in terms of the deformation theory of free \textit{c.d.g.a's
}in characteristic zero. The discussion that follows makes the connection precise.

Throughout the remainder of this section, $\mathbf{k}$ denotes a field of
characteristic zero. Let $(A,d,\mu)$ be a free $c.d.g.a.$ over $\mathbf{k}$
and let $(M,\lambda,\rho,d_{M})$ be a symmetric \textit{d.g. }$A$-bimodule. By
forgetting the differentials, Theorem \ref{free} implies that $Harr^{n}%
(A,M)=0$ for $n>1$ so that each column in the \textit{restricted} Harrison
complex is exact (see Figure 3).\bigskip\ %

\[%
\begin{array}
[c]{c}%
\begin{array}
[c]{ccccc}
&  & \;\vdots &  & \\
&  & {\scriptstyle\partial}_{\scriptscriptstyle B}\uparrow\quad &  & \\
&  &  &  & \\
\cdots & \rightarrow & Ch^{p,3}(A;M) & \rightarrow & \cdots\\
&  &  &  & \\
&  & {\scriptstyle\partial}_{\scriptscriptstyle B}\uparrow\quad &  & \\
&  &  &  & \\
\cdots & \rightarrow & Ch^{p,2}(A;M) & \rightarrow & \cdots\\
&  &  &  & \\
&  & {\scriptstyle\partial}_{\scriptscriptstyle B}\uparrow\quad &  & \\
&  &  &  & \\
\cdots & \rightarrow & Ch^{p,1}(A;M) & \rightarrow & \cdots
\end{array}
\\
\\
\text{An Exact Column}\\
\text{Figure 3.}%
\end{array}
\]
\vspace{0.2in}

Let $\{E,d_{E}^{\prime},d_{E}^{\prime\prime}\}$ be a standard double cochain
complex with differentials of respective bidegree $(1,0)$ and $(0,1).$

\begin{definition}
Let $n\geq1,$ $n\geq k\geq0,$ and let $f=\sum_{0\leq i\leq n-1}f_{i,n-i}$ be a
total $n$-cocycle for which $\{f_{i,n-i}\in E^{i,n-i}\}_{0\leq i\leq n-1}.$
Then $f$ is an $(n,k)$-cocycle if $f_{i,n-i}=0$ for $0\leq i\leq k-1.$ An
$(n,n)$-cocycle is said to be concentrated in bidegree $(n-1,1).$
\end{definition}

\begin{lemma}
\label{lemma}If the columns of $\{E,d_{E}^{\prime},d_{E}^{\prime\prime}\}$ are
exact, then every class $[x]\in H^{n}(E,d_{E}^{\prime}+d_{E}^{\prime\prime}) $
can be represented by a total $n$-cocycle $h$ concentrated in bidegree
$(n-1,1).$
\end{lemma}

\textit{Proof:} Let $f=\sum_{0\leq i\leq n-1}f_{i,n-i}$ be a total $n$-cocycle
and note that $f$ is an $(n,k)$-cocycle for some $k=0,1,2,\ldots n.$ If $k=n$
there is nothing to prove; so assume that $k<n.$ By exactness, there exists a
cochain $g_{k,n-k-1}\in E^{k,n-k-1}$ such that $d_{E}^{\prime\prime
}(g_{k,n-k-1})=f_{k,n-k}$ and $d_{E}^{\prime\prime}$$(f_{k+1,n-k-1}%
-d_{E}^{\prime}(g_{k,n-k-1}))=0.$ Hence $h_{k+1}=\newline\sum_{k+1\leq i\leq
n-1}f_{i,n-i}-d_{E}^{\prime}(g_{k,n-k-1})$ is an $(n,k+1)$-cocycle and
$f-h_{k+1}=$\linebreak$d_{E}^{\prime}(g_{k,n-k-1})+f_{k,n-k}=(d_{E}^{\prime
}+d_{E}^{\prime\prime})(g_{k,n-k-1})$. Proceed inductively until $f$ is
totally cohomologous to some $(n,n)$-cocycle $h_{n}$.\bigskip\ 

I should note that ''staircase'' arguments such as this are not new, having
appeared as early as 1952 in a paper by Weil \cite{Weil}.

\begin{definition}
Let $(M,\lambda,\rho)$ be an $A$-bimodule. A $\mathbf{k}$-linear map
$\theta:A\rightarrow M$ is a derivation if $\theta\circ\mu=\rho\circ
(\theta\otimes1)+\lambda\circ(1\otimes\theta).$ The set of all derivations of
degree $p$ is denoted by $Der^{p}(A,M).$
\end{definition}

\begin{corollary}
\label{der}Let $(A,d,\mu)$ be a free $c.d.g.a.$ over $\mathbf{k}$ and let
$(M,\lambda,\rho,d_{M})$ be a symmetric d.g. $A$-bimodule. If $f=\sum_{0\leq
i\leq n-1}f_{i,n-i}\in\widetilde{C}h^{n}(A;M;3)$ is a total $n$-cocycle, there
exists $g\in Der^{n-1}(A,M)$ such that $f-g$ is totally cohomologous to zero.
\end{corollary}

\textit{Proof:} By Lemma \ref{lemma}, there exists $n$-cocycle $g$
concentrated in bidegree $(n-1,1)$ such that $f-g$ is totally cohomologous to
zero$.$ But then $g\in ker$\text{\textit{\ }}$\partial_{B}$ and $Der^{n-1}%
(A,M)=ker$\text{\textit{\ }}$\partial_{B}$ by the definition of $\partial
_{B}.$

\begin{theorem}
Let $(A,d,\mu)$ be a free $c.d.g.a.$ over $\mathbf{k}$ and let $(M,\lambda
,\rho,d_{M})$ be a symmetric \textit{d.g. }$A$-bimodule. There is an
isomorphism
\begin{equation}
\Gamma:H^{*-1}(Der(A,M),d_{B})\overset{\approx}{\rightarrow}\widetilde
{H}arr_{c.d.g.a.}^{*}(A;M;3).\label{iso}%
\end{equation}

\end{theorem}

\textit{Proof:} Consider a $d_{B}$-cocycle $\theta\in Der^{n-1}(A,M);$
$\theta$ is a normalized $(n-1)$-cochain since $\theta(1_{\mathbf{k}})=0.$ On
the other hand, $\theta$ is an $(n,n)$-cocycle in $\widetilde{C}h^{n}(A,M;3).$
So define $\Gamma[\theta]=[\theta];$ surjectivity follows from Corollary
\ref{der}.\bigskip\ 

Since every Harrison $n$-class $[x]$ can be represented by a total $n$-cocycle
concentrated in bidegree $(n-1,1),$ the dimension shift in isomorphism
(\ref{iso}) is superficial to the extent that it emphasizes one of two
points-of-view. We can think of $[x]$ either in terms of a representative in
$Der^{n-1}(A;M)\subset Ch^{n-1,1}(A;M)$ or in terms of a representative in
$\widetilde{C}h^{n}(A;M;3).$

Now set $M=A.$ The adjoint action of $d$ on $Der^{*}(A,A)$ as a Lie algebra of
derivations is given by $ad(d)(\theta)=[d,\theta],$ where $[d,\theta
]=d\circ\theta-(-1)^{|\theta|}\theta\circ d$. In this case, the differential
in (\ref{iso}) is simply
\[
d_{B}=-ad(d).
\]
In particular, let $\mathbf{\Lambda}(x_{i})$ denote the free \textit{c.g.a.
}on generators $\{x_{i}\}$ over $\mathbf{k}$ and let $\left\langle
x_{i}\right\rangle $ denote the $\mathbf{k}$-module with basis $\{x_{i}\}$.
Since derivations of a free \textit{c.g.a}. are determined by their action on
generators and $\mathbf{k}$-linear maps are determined by their action on a
basis we have:

\begin{corollary}
If $d$ is an algebra differential on $A=\mathbf{\Lambda}(x_{i})$, there is an
isomorphism
\begin{equation}
\widetilde{H}arr_{c.d.g.a.}^{*}(A;A;3)\approx H(Hom_{\mathbf{k}}%
^{*-1}(\left\langle x_{i}\right\rangle ,A),ad(d)).\label{iso2}%
\end{equation}

\end{corollary}

Note that the right-hand-side of isomorphism (\ref{iso2}) does not depend upon
the multiplication $\mu.$ Thus a non-vanishing Harrison 2-class $[x]\in
\widetilde{H}arr_{c.d.g.a.}^{2}(A;A;3)$\textit{\ }signals a potential change
in the differential $d$---\textit{not }in the multiplication $\mu$. But
changing the differential is exactly the game Felix \cite{Felix} and
Halperin-Stasheff \cite{Halperin} play.

Let $X$ be a formal space with rational cohomology algebra $A=H^{*}%
(X;\mathbb{Q}\mathbf{)}.$ A \textit{rational minimal model for }$A$ is a free
\textit{c.d.g.a. }$(\Lambda,d)$ over $\mathbb{Q}$ such that $d(\Lambda
)\subset\Lambda\cdot\Lambda$ and $H^{*}(\Lambda,d)\approx A.$ For Felix, this
isomorphism is additive; for Halperin and Stasheff it is multiplicative. A
\textit{perturbation} $p$ of the differential $d$ is a derivation $p\in
Der^{1}(\Lambda,\Lambda)$ such that $(d+p)^{2}=0$ and $H^{*}(\Lambda
,d+p)\approx A$. The linearization of a perturbation $p$ represents a
\textit{1-class} on the right-hand-side of (\ref{iso2}). Given a perturbation
$p,$ there exists a rational space $Y$ and an isomorphism of higher order
structures $H^{*}(\Lambda,d+p)\cong H^{*}(Y;\mathbb{Q})$. Conversely, given a
rational space $Y$ with $H^{*}(Y;\mathbb{Q})\approx A,$ there exists a
perturbation $p_{Y}$ such that $H^{*}(\Lambda,d+p_{Y})\cong H^{*}%
(Y;\mathbb{Q})$. Spaces $Y$ and $Y^{\prime}$ have the same rational homotopy
type if and only if the corresponding perturbations $p_{Y}$ and $p_{Y^{\prime
}}$ are equivalent as deformations (see section 6 below).

Felix and Halperin-Stasheff apply this theory in two somewhat different ways.
One can obtain a minimal model for $A$ as the limit of a Tate-Jozefiak
resolution of $A$ \cite{Tate}, \cite{Jozefiak}. Thought of this way, the
minimal model is naturally bigraded with respect to internal and resolution
degrees. In either approach, this ''bigraded model'' is perturbed to a
''filtered model''; the perturbations of Felix arbitrarily decrease
filtration, while those of Halperin and Stasheff decrease filtration by at
least two. The former effectively fixes the additive structure and varies the
multiplicative structure on $A$, while the latter fixes the multiplicative
structure and varies the higher order algebra structure on $A$. The set of
perturbations that decrease filtration by at least two is a subcomplex of
$\{Der^{*}(\Lambda,\Lambda),$ $ad(d)\}$.

Finally, I should mention that a comparison between the Harrison cohomology of
$(\Lambda,d)$ and the homology of $\{Coder(L^{c}\Lambda,A),d^{c}\},$ where
$L^{c}\Lambda$ is the free \textit{d.g. }Lie coalgebra on $\Lambda$ and
$d^{c}$ is the differential induced by $d,$ was given by Schlessinger and
Stasheff in \cite{Schlessinger}.\bigskip

\subsection{The Deformation Complex for Differential Graded Coalgebras}

Let $(C,d,\Delta)$ be a \textit{d.g.c. }The deformation complex for $C$ is a
double complex dual to the one discussed in 5.1; this discussion is included
for notational purposes. Let $(N,\lambda,\rho,d_{N})$ be a \textit{d.g.}
$C$-bicomodule$.$ For each $n\geq0,$ consider the exterior $C$-bicomodule
$(C^{\otimes(n+2)},d_{(n)})$ where $d_{(n)}$ is defined as in (\ref{d-sub-n}).
By proposition \ref{bicomod-map} and an easy calculation, there is a map
$d^{q,n}:Hom_{C-bicomod}^{q}(N,C^{\otimes(n+2)})\rightarrow\newline
Hom_{C-bicomod}^{q+1}(N,C^{\otimes(n+2)})$ given by
\[
d^{q,n}(g)=(-1)^{q}g\circ d_{N}-d_{(n)}\circ g\text{\ .}
\]
This induces a $\mathbf{k}$-linear differential $d_{\Omega}^{q,n}%
:Hom_{\mathbf{k}}^{q}(N,C^{\otimes n})\rightarrow Hom_{\mathbf{k}}%
^{q+1}(N,C^{\otimes n})$ of bidegree $(1,0)$ via (\ref{psi}):
\[
d_{\Omega}^{*,*}=\text{ }\Psi\circ d^{*,*}\circ\Psi^{-1}\text{.}
\]
Via (\ref{XI-star}), the $C$-bicomodule map $\delta_{(n)}:C^{\otimes
(n+2)}\rightarrow C^{\otimes(n+3)}$ induces a $\mathbf{k}$-linear differential
$\delta_{\Omega}^{q,n}:Hom_{\mathbf{k}}^{q}(N,C^{\otimes n})\rightarrow
Hom_{\mathbf{k}}^{q}(N,C^{\otimes(n+1)})$ of bidegree $(1,0)$ given by
\[
\delta_{\Omega}^{*,*}=\Psi\circ\delta^{*,*}\circ\Psi^{-1}.
\]
It is easy to check that $d_{\Omega}\circ\delta_{\Omega}-\delta_{\Omega}\circ
d_{\Omega}=0;$ at a particular $g\in Hom_{\mathbf{k}}^{q}(N,C^{\otimes n})$ we
have
\begin{equation}
d_{\Omega}^{q,n}(g)=(-1)^{q}g\circ d_{N}-d_{(n-2)}\circ g\text{ }\label{dc}%
\end{equation}
and
\begin{equation}
\delta_{\Omega}^{q,n}(g)=(1\otimes g)\circ\lambda-\delta_{(n-2)}\circ
g+(-1)^{n+1}(g\otimes1)\circ\rho,\label{d-cobar}%
\end{equation}
where $d_{(-2)}=0$ and $\delta_{(-2)}=0.$

Refer to $\Omega^{q,n}(C;N)=Hom_{\mathbf{k}}^{q}(N,C^{\otimes n})$ as the
space of \textit{Cartier n-cochains on }$C$ \textit{of degree q. }The double
complex $\{\Omega^{*,*}(C;N),d_{\Omega},\delta_{\Omega}\}$ is called the
\textit{Cartier cochain complex on }$C$\textit{\ with coefficients in the d.g.
A-bicomodule }$N$. Define the space of total $s$-cochains by $\Omega
^{s}(C;N)=\sum_{q\in\mathbb{Z}\mathbf{;}\text{ }s\geq q}\Omega^{q,s-q}%
(C;N)\ $and define $D_{\Omega}$ on the component $\Omega^{q,n}(C;N)$ by
\[
D_{\Omega}=\delta_{\Omega}-(-1)^{n}d_{\Omega}.
\]
Then $D_{\Omega}$ is a total differential; the sign $-(-1)^{n}$ is introduced
so that $D_{\Omega}^{2}=0.$ Now define the \textit{Cartier cohomology of }$C$
\textit{with coefficients in the d.g. A-bicomodule }$N,$ denoted by
$H_{d.g.c}^{*}(C;N),$ to be the homology of the total complex $\{\Omega
^{*}(C;N),D_{\Omega}\}$.

We say that an $n$-cochain $g\in\Omega^{q,n}(C;N)$ is \textit{normalized} if
$g(x)=a_{1}\otimes\cdots\otimes a_{n}=0$ whenever $\left|  a_{i}\right|  =0$
for some $i=1,2,\ldots n.$ The differential $\delta_{\Omega}$ restricts to the
subspace of normalized cochains and a standard theorem \cite{MacLane} assures
that the subcomplex given by such a restriction is cochain homotopic to the
full cochain complex. Consequently, we shall use normalized cochains; the
symbol $\Omega^{q,n}(C;N)$ will henceforth denote the space of normalized
$n$-cochains of degree $q$.

As on the algebra side, obtain the \textit{restricted }Cartier cohomology of
$C$ by deleting the bottom row of the Cartier cochain complex$.$ Denote the
space of restricted total cochains by $\widetilde{\Omega}^{*}(C;N)$ and the
restricted cohomology by $\widetilde{H}_{d.g.c.}^{*}(C;N).$ For $m\in
\{3,4,\ldots\}\cup\{\infty\},$ further restrict the Cartier complex to those
cochains in bidegree $(q,n)$ with $q\geq3-m$ and $n\geq1.$ Denote the
bitruncated total $s$-cochains by $\widetilde{\Omega}^{s}(C;N;m)=\sum_{s-1\geq
q\geq3-m}\Omega^{q,s-q}(C;N).$ The homology of the complex $\{\widetilde
{\Omega}^{*}(C;N;m),\newline D_{\Omega}\},$ denoted by $\widetilde{H}%
_{d.g.c.}^{*}(C;N;m),$ is called the \textit{restricted Cartier cohomology of
}$C$\textit{\ truncated at degree }$3-m$\textit{. }

The \textit{deformation complex for $C$ as a strict }$C\left(  m\right)
$\textit{-coalgebra }is the cochain complex $\{\widetilde{\Omega}^{\ast
}(C;C;m),D_{\Omega}\}.$ The cohomology $\widetilde{H}_{d.g.c.}^{\ast}(C;C;m)$
directs the deformation theory in a way completely analogous to the algebra
case. If $C_{t}=(C[[t]],\ \Delta_{t}^{(1)},\ \Delta_{t}^{(2)},$\linebreak%
$\ldots)$ is a deformation of $C$ as a strict $C(m)$-coalgebra, we agree that
for $1\leq i\leq m$ the maps $\Delta_{t}^{(i)}=\Delta_{0}^{(i)}+t\Delta
_{1}^{(i)}+t^{2}\Delta_{2}^{(i)}+\cdots$ satisfy $\Delta_{1}^{(i)}\in
\Omega^{2-i,i}(C;C)$ with $\Delta_{0}^{(1)}=d,$ $\Delta_{0}^{(2)}=\Delta,$ and
$\Delta_{0}^{(i)}=0$ for $i>2.$ Then $c=\sum_{i=1}^{m}\Delta_{1}^{(i)}\in$
$\widetilde{\Omega}^{2}(C;C;m)$ satisfies $D_{\Omega}(c)=0.$ Furthermore,
given a total cocycle $c=\sum_{i=1}^{m}\Delta_{1}^{(i)}\in$ $\widetilde
{\Omega}^{2}(C;C;m)$, the obstructions to extending the corresponding linear
approximation $(d+t\Delta_{1}^{(1)},$ $\mu+t\Delta_{1}^{(2)},$ $\ldots
,\ t\Delta_{1}^{(i)},\ldots)$ to a deformation appear as an inductively
defined sequence cocycles in $\widetilde{\Omega}^{3}(C;C;m).$ In the special
case $m=3$ we adopt the standard notation $C_{t}=(C[[t]],\ d_{t},\ \Delta
_{t})$ with $d_{t}=d+td_{1}+t^{2}d_{2}+\cdots$ and $\Delta_{t}=\Delta
+t\Delta_{1}+t^{2}\Delta_{2}+\cdots,$ in which case $C_{t}$ is a
(coassociative) \textit{d.g.c. }We refer further discussion of this latter
case to section 6.

\subsection{The Deformation Complex for Differential Graded Hopf Algebras}

Let $(H,d,\mu,\Delta)$ be a \textit{d.g.h.a.} In Example \ref{d.g.bidimod} we
observed that $(H^{\underline{\otimes}(m+2)},d_{(m)})$ and $(H^{\overline
{\otimes}(n+2)},d_{(n)})$ are \textit{d.g}. $H$-bidimodules and that
$\partial_{(m)}:H^{\underline{\otimes}(m+2)}\rightarrow H^{\underline{\otimes
}(m+1)}$ and $\delta_{(n)}:H^{\overline{\otimes}(n+2)}\rightarrow
H^{\overline{\otimes}(n+3)}$ are \textit{d.g. }$H$-bidimodule maps. By
proposition \ref{bidimod-map}, there is a map $d^{p,m,n}:Hom_{H-bidimod}%
^{p}(H^{\underline{\otimes}(m+2)},H^{\overline{\otimes}(n+2)})\rightarrow
$\linebreak$Hom_{H-bidimod}^{p+1}(H^{\underline{\otimes}(m+2)},H^{\overline
{\otimes}(n+2)})$ given by $d^{p,m,n}(f)=(-1)^{p}f\circ d_{(m)}-d_{(n)}\circ
f.$

In section 5.1 and 5.2 we observed that $d^{*,*,*}$ commutes with the induced
maps $\partial^{p,m,n}=Hom_{H-bimod}^{p}(\partial_{(m)},H^{\overline{\otimes
}(n+2)})$ and $\delta^{p,m,n}=Hom_{H-bicomod}^{p}(H^{\underline{\otimes}%
(m+2)},\delta_{(n)});$ consequently $d^{*,*,*}$ commutes with $\partial
^{*,*,*}$ and $\delta^{*,*,*}$ in the common subcategory of $H$-bidimodules.
But there is more. For each $p\in\mathbb{Z}$ and each $m,n\geq0,$ the bar and
cobar differentials (\ref{bar}) and (\ref{cobar}) functorially induce the
commutativity of $\partial^{*,*,*}$ and $\delta^{*,*,*}$ in this subcategory
(see Figure 4).This remarkable structural compatibility gives rise to a triple
complex $\{Hom_{H-bidimod}^{*}(H^{\underline{\otimes}(*+2)},H^{\overline
{\otimes}(*+2)}),$ $d^{*,*,*},$ $\partial^{*,*,*},$ $\delta^{*,*,*}\}.$ The
desired triple complex at the level of $\mathbf{k}$-modules is obtained via
isomorphism (\ref{isobidimod}) with differentials induced by $d^{*,*,*},$
$\partial^{*,*,*},$ and $\delta^{*,*,*}$as in (\ref{THETA-STAR}) and
(\ref{XI-STAR}). We have obtained our main result.

\begin{theorem}
Let $(H,d,\mu,\Delta)$ be a \textit{d.g.h.a. }and let $C^{p,m,n}%
(H;H)=$\linebreak$Hom_{\mathbf{k}}^{p}(H^{\underline{\otimes}m},H^{\overline
{\otimes}n})$. There is a triple complex $\{C^{\ast,\ast,\ast}(H;H),d_{C}%
,\partial_{C},\delta_{C}\}$ whose differentials have respective tridegree
$(1,0,0),$ $(0,1,0),$ and $(0,0,1)$ and arise from $d,$ $\mu,$ and $\Delta$ as
follows:
\[
d_{C}^{\ast,\ast,\ast}=(\Phi\circ\Psi)\circ d^{\ast,\ast,\ast}\circ(\Phi
\circ\Psi)^{-1},
\]%
\[
\partial_{C}^{\ast,\ast,\ast}=(\Phi\circ\Psi)\circ\partial^{\ast,\ast+1,\ast
}\circ(\Phi\circ\Psi)^{-1},\text{ and}%
\]%
\[
\delta_{C}^{\ast,\ast,\ast}=(\Phi\circ\Psi)\circ\delta^{\ast,\ast,\ast}%
\circ(\Phi\circ\Psi)^{-1}.
\]
At a particular $f\in C^{p,m,n}(H;H)$ these expand to
\[
d_{C}^{p,m,n}(f)=(-1)^{p}f\circ d_{(m-2)}-d_{(n-2)}\circ f,
\]%
\[
\partial_{C}^{p,m,n}{}(f)=\lambda^{n}\circ(1\otimes f)-f\circ\partial
_{(m-1)}+(-1)^{m+1}\rho^{n}\circ(f\otimes1),\text{ \ and}%
\]%
\[
\delta_{C}^{p,m,n}(f)=(1\otimes f)\circ\lambda_{m}-\delta_{(n-2)}\circ
f+(-1)^{n+1}(f\otimes1)\circ\rho_{m}.
\]

\end{theorem}

\noindent$%
\begin{array}
[c]{c}%
\begin{array}
[c]{ccc}
&  & \\
& \partial^{p,m,n+1} & \\
Hom_{H-bidimod}^{p}(H^{\underline{\otimes}m+2},H^{\overline{\otimes}n+3}) &
\longrightarrow & Hom_{H-bidimod}^{p}(H^{\underline{\otimes}m+3}%
,H^{\overline{\otimes}n+3})\\
&  & \\
\delta^{p,m,n}\uparrow &  & \uparrow\delta^{p,m+1,n}\\
&  & \\
Hom_{H-bidimod}^{p}(H^{\underline{\otimes}m+2},H^{\overline{\otimes}n+2}) &
\longrightarrow & Hom_{H-bidimod}^{p}(H^{\underline{\otimes}m+3}%
,H^{\overline{\otimes}n+2})\\
& \partial^{p,m,n} &
\end{array}
\\
\\
\text{Figure 4.}%
\end{array}
\vspace{0.2in}$

The space $C^{p,m,n}(H;H)$ of normalized cochains is referred to as the
\textit{Hochschild-Cartier (m,n)-cochains on }$H$ \textit{of degree p; }the
triple complex $\{C^{\ast,\ast,\ast}(H;H),d_{C},\partial_{C},$\linebreak%
$\delta_{C}\}$ is called the \textit{Hochschild-Cartier cochain complex on
}$H.$ Define the space of total $r$-cochains by $C^{r}(H;H)=\sum
_{p\in\mathbb{Z};\text{ }m,n\geq0;\text{ }p+m+n=r+1}C^{p,m,n}(H;H)$ and define
$D_{C}$ on the component $C^{p,m,n}(H;H)$ by
\[
D_{C}=(-1)^{m(n+1)}\;d_{C}+(-1)^{n(p+1)}\partial_{C}+(-1)^{p(m+1)}\delta_{C}.
\]
Then $D_{C}$ is a total differential; the sign adjustments are introduced so
that $D_{C}^{2}=0$ and determine those in (\ref{d}), (\ref{d-bar}),
(\ref{dc}), and (\ref{d-cobar}) made earlier---take coefficients in $H$ and
set either\textit{\ }$m=1$ or $n=1$ as appropriate. Finally, define the
\textit{Hochschild-Cartier cohomology of }$H$\textit{, denoted by
}$H_{d.g.h.a.}^{\ast}(H;H),$ to be the homology of the complex $\{C^{\ast
}(H;H),D_{C}\}.$

Obtain the \textit{restricted }Hochschild-Cartier cohomology of $H$ by
deleting the two ''coordinate planes'' $m=0$ and $n=0$ in the
Hochschild-Cartier triple cochain complex. Denote the space of restricted
cochains by $\widetilde{C}^{*}(H;H);$ denote the corresponding restricted
cohomology by $\widetilde{H}_{d.g.h.a.}^{*}(H;H).$ For $q\in\{3,4,\ldots
\}\cup\{\infty\},$ further restrict the Hochschild-Cartier complex to those
cochains in tridegree $(p,m,n)$ with $p\geq3-q$ and $m,n\geq1.$ Denote the
space of tritruncated total $r$-cochains by $\widetilde{C}^{r}(H;H;q)=\sum
_{m,n\geq1;\text{ }p\geq3-q;\text{ }p+m+n=r+1}C^{p,m,n}(H;H).$

The \textit{deformation complex for }$H$ \textit{as a strict }$H(q)$%
-\textit{structure} is the complex\newline$\{\widetilde{C}^{\ast}%
(H;H;q),D_{C}\}$\textit{;} its homology, which directs the deformation theory,
is called the \textit{restricted Hochschild-Cartier cohomology of H truncated
at degree }$3-q$ and is denoted by $\widetilde{H}_{d.g.h.a.}^{\ast}(H;H;q)$.
When $q=3$ we obtain the deformation complex for $H$ as a \textit{d.g.h.a. }If
the Hopf algebra $H$ is commutative as an algebra, we obtain the
\textit{restricted Harrison cohomology of }$H$ \textit{truncated at degree}
$3-q,$ denoted by $\widetilde{H}arr_{c.d.g.h.a.}^{\ast}(H;H;q),$ as the
homology of $\{\widetilde{C}h^{\ast}(H;H;q),D_{C}\}$. We conclude our
discussion with a brief exposition of the deformation theory of
\textit{d.g.h.a}'s.\ 

\section{Deformation Theory of Differential Graded Hopf Algebras}

Let $H_{0}=(H,d,\mu,\Delta)$ be a \textit{d.g.h.a}. over a field $\mathbf{k}.$
Let $t$ be an indeterminant of degree $0$ and let $\mathbf{k}\left[  \left[
t\right]  \right]  $ denote the commutative ring of formal power series in
$t.$ Consider the (graded) $\mathbf{k}\left[  \left[  t\right]  \right]
$-module $H\left[  \left[  t\right]  \right]  $ of formal power series in $t$
with coefficients in $H.$ We give $H[[t]]$ the $t$-\textit{adic topology }in
which $a $ and $b$ are $t^{r}$\textit{-close} if $a\equiv b(mod\ t^{r}).$ In
the $t$-adic topology, every formal power series is the limit of its sequence
of partial sums.

Given a $\mathbf{k}$-linear map $f:H^{\underline{\otimes}m}\rightarrow
H^{\overline{\otimes}n}$, extend $f$ to a $\mathbf{k}[[t]]$-linear map
$f:H[[t]]^{\underline{\otimes}m}\rightarrow H[[t]]^{\overline{\otimes}n},$
where we tensor over $\mathbf{k}[[t]],$ by defining $f(\sum t^{i}a_{i}%
\otimes\sum t^{j}a_{j}\otimes\cdots)=\sum t^{i+j+\cdots}f(a_{i}\otimes
a_{j}\otimes\cdots);$ this is the unique $\mathbf{k}\left[  \left[  t\right]
\right]  $-linear extension of $f$ to $H[[t]]^{\underline{\otimes}m}.$ In
particular, so extending $d,$ $\mu$ and $\Delta$ gives a $\mathbf{k}\left[
\left[  t\right]  \right]  $-\textit{d.g.h.a} $H_{0}[[t]]=(H\left[  \left[
t\right]  \right]  ,d,\mu,\Delta),$ which possibly deforms to some
$\mathbf{k}\left[  \left[  t\right]  \right]  $-\textit{d.g.h.a}.
$H_{t}=(H\left[  \left[  t\right]  \right]  ,d_{t},\mu_{t},\Delta_{t})$.
Indeed, the deformation theories of algebras and coalgebras discussed above
suggest the possibility of further deforming to some $H_{t}$ in the category
of ''$H(q)$-bialgebras'', but we limit our discussion to the deformation of
$H_{0}$ as a \textit{d.g.h.a. }here.

\begin{definition}
A deformation of $H_{0}$ is a $\mathbf{k}\left[  \left[  t\right]  \right]
$-\textit{d.g.h.a}. $H_{t}=(H\left[  \left[  t\right]  \right]  ,d_{t},\mu
_{t},$\linebreak$\Delta_{t})$ such that

\begin{enumerate}
\item $d_{t}=d+td_{1}+t^{2}d_{2}+\cdots,$

\item $\mu_{t}=\mu+t\mu_{1}+t^{2}\mu_{2}+\cdots,$ and

\item $\Delta$$_{t}=\Delta+t\Delta_{1}+t^{2}\Delta_{2}+\cdots$.
\end{enumerate}
\end{definition}

The fundamental problem is to classify all \textquotedblright
non-trivial\textquotedblright\ deformations. Given a deformation
$H_{t}=(H\left[  \left[  t\right]  \right]  ,d_{t},\mu_{t},\Delta_{t}),$
consider the linear terms in the associativity condition $\mu_{t}%
\circ(1\otimes\mu_{t})=\mu_{t}\circ(\mu_{t}\otimes1).$ Equating coefficients
gives $\mu\circ(1\otimes\mu_{1})+$\linebreak$\mu_{1}\circ(1\otimes\mu
)=\mu\circ(\mu_{1}\otimes1)+\mu_{1}\circ(\mu\otimes1)$ $.$ With $\lambda
^{1}=\rho^{1}=\mu$ we have $\partial_{C}^{0,2,1}(\mu_{1})=\lambda^{1}%
\circ(1\otimes\mu_{1})-\mu_{1}\circ(\mu\otimes1)+\mu_{1}\circ(1\otimes
\mu)-\rho^{1}\circ(\mu_{1}\otimes1)=0.$ Similarly, with $\lambda_{1}=\rho
_{1}=\Delta,$ the linear terms in the coassociativity condition $(1\otimes
\Delta_{t})\circ\Delta_{t}=(\Delta_{t}\otimes1)\circ\Delta_{t}$ imply that
$\delta_{C}^{0,1,2}$$(\Delta_{1})=0,$ and the linear terms in the differential
condition $d_{t}\circ d_{t}=0$ imply that $d_{C}^{1,1,1}(d_{1})=0.$
Furthermore, the linear terms in the bialgebra condition $\Delta_{t}\circ
\mu_{t}=(\mu_{t}\otimes\mu_{t})\circ(2,3)\circ(\Delta_{t}\otimes\Delta_{t})$
imply that $\partial_{C}^{0,1,2}$$(\Delta_{1})+\delta_{C}^{0,2,1}(\mu
_{1})=\lambda^{2}\circ(1\otimes\Delta_{1})-\Delta_{1}\circ\mu+\rho^{2}%
\circ(\Delta_{1}\otimes1)+(1\otimes\mu_{1})\circ\lambda_{2}-\Delta\circ\mu
_{1}+(\mu_{1}\otimes1)\circ\rho_{2}$ $=$ $0;$ the linear terms in the
derivation condition $d_{t}\circ\mu_{t}=\mu_{t}\circ(d_{t}\otimes1+1\otimes
d_{t})$ imply that $d_{C}^{0,2,1}(\mu_{1})+\partial_{C}^{1,1,1}(d_{1})=0;$ and
the linear terms in the coderivation condition $\Delta_{t}\circ d_{t}%
=(d_{t}\otimes1+1\otimes d_{t})\circ\Delta_{t}$ imply that $d_{C}%
^{0,1,2}(\Delta_{1})-\delta_{C}^{1,1,1}(d_{1})=0.$ Hence $D_{C}(d_{1}+\mu
_{1}+\Delta_{1})=0$ so that $d_{1}+\mu_{1}+\Delta_{1}$ is a restricted total
Hochschild-Cartier $2$-cocycle. Restricting these calculations to
\textquotedblright planes\textquotedblright\ $n=1$ and $m=1$ establishes
similar claims made earlier in the algebra and coalgebra settings, respectively.

Conversely, given a restricted total Hochschild-Cartier $2$-cocycle $d_{1}%
+\mu$$_{1}+\Delta_{1}\in\newline C^{1,1,1}(H;H)\oplus C^{0,2,1}(H;H)\oplus
C^{0,1,2}(H;H), $ there is an inductively defined sequence of cocycles in
$\widetilde{C}^{3}(H;H;3)$ whose vanishing allows one to inductively extend
the corresponding linear approximation $(d+td_{1},$ $\mu+t\mu_{1},$
$\Delta+t\Delta_{1})$ to a deformation $H_{t}$ in a standard way. For a
detailed discussion of the bialgebra case see \cite{Gersten}.

Recall that a map $\phi_{t}\in Aut_{\mathbf{k}[[t]]}(H[[t]])\ $satisfies
$\phi_{t}\equiv1(mod\ t);$ hence there exist $\mathbf{k}$-linear maps
$\{\phi_{i}\}_{i\geq1}$ such that $\phi_{t}=1+t\phi_{1}+t^{2}\phi_{2}+\cdots.$

\begin{definition}
\label{equival}Two deformations $H_{t}=(H\left[  \left[  t\right]  \right]
,d_{t},\mu_{t},\Delta_{t})$ and $H_{t}^{\prime}=(H\left[  \left[  t\right]
\right]  ,$\linebreak$d_{t}^{\prime},\mu_{t}^{\prime},\Delta_{t}^{\prime})$
are equivalent\textit{\ }if there exists some \textit{d.g.h.a.} map $\phi
_{t}\in Aut_{\mathbf{k}[[t]]}(H[[t]])$ such that

\begin{enumerate}
\item $d_{t}\circ\phi_{t}=\phi_{t}\circ d_{t}^{\prime},$

\item $\mu_{t}\circ(\phi_{t}\otimes\phi_{t})=\phi_{t}\circ\mu_{t}^{\prime},$ and

\item $\Delta_{t}\circ\phi_{t}=(\phi_{t}\otimes\phi_{t})\circ\Delta
_{t}^{\prime}.$
\end{enumerate}

\noindent When $H_{t}$ and $H_{t}^{\prime}$ are equivalent via $\phi_{t}$, we
write $\phi_{t}:H_{t}\sim H_{t}^{\prime}$ and refer to $\phi_{t}$ as an
equivalence. Any deformation $H_{t}$ equivalent to $H_{0}[[t]]=(H[[t]],d,\mu
,\Delta),$ thought of as a deformation with $d_{i}=\mu_{i}=\Delta_{i}=0$ for
all $i,$ is called a trivial deformation. If every deformation $H_{t}$ is
trivial we say that $H_{0}$ is rigid as a d.g.h.a.
\end{definition}

Given an equivalence $\phi_{t}:H_{t}\sim H_{t}^{\prime}$, equate coefficients
in the linear terms of the naturality condition to obtain $d_{C}^{0,1,1}%
(\phi_{1})=\phi_{1}\circ d-d\circ\phi_{1}=d_{1}-d_{1}^{\prime}.$ Similarly,
$-\partial_{C}^{0,1,1}(\phi_{1})=-[\mu\circ(1\otimes\phi_{1})-\phi_{1}\circ
\mu+\mu\circ(\phi_{1}\otimes1)]=\mu_{1}-\mu_{1}^{\prime}$ since $\phi_{t}$ is
an algebra map, and $\delta_{C}^{0,1,1}(\phi_{1})=(1\otimes\phi_{1}%
)\circ\Delta-\Delta\circ\phi_{1}+(\phi_{1}\otimes1)\circ\Delta=\Delta
_{1}-\Delta_{1}^{\prime}$ since $\phi_{t}$ is a coalgebra map$.$ Therefore
$D_{C}(\phi_{1})=(d_{1}+\mu_{1}+\Delta_{1})-(d_{1}^{\prime}+\mu_{1}^{\prime
}+\Delta_{1}^{\prime})$ so that $d_{1}+\mu_{1}+\Delta_{1}$ and $d_{1}^{\prime
}+\mu_{1}^{\prime}+\Delta_{1}^{\prime}$ are totally cohomologous as restricted
Hochschild-Cartier $2$-cocycles.

If $\widetilde{H}_{d.g.h.a.}^{2}(H;H;3)=0$ and $H_{t}=(H[[t]],d_{t},\mu
_{t},\Delta_{t})$ is a deformation, choose $\phi_{1}\in\widetilde{C}%
^{1}(H;H;3)$ such that $D_{C}^{0,1,1}(\phi_{1})=d_{1}+\mu_{1}+\Delta_{1}$ and
consider $\phi_{t}^{(1)}=(1-t\phi_{1})\in Aut_{\mathbf{k}[[t]]}(H[[t]]).$
There is a deformation $H_{t}^{(1)}=(H\left[  \left[  t\right]  \right]
,d_{t}^{(1)},\mu_{t}^{(1)},\Delta_{t}^{(1)}),$ where $d_{t}^{(1)}=\phi
_{t}^{(1)}\circ d_{t}\circ[\phi_{t}^{(1)}]^{-1},$ $\mu_{t}^{(1)}=\phi
_{t}^{(1)}\circ d_{t}\circ[\phi_{t}^{(1)}\otimes\phi_{t}^{(1)}]^{-1},$ and
$\Delta_{t}^{(1)}=[\phi_{t}^{(1)}\otimes\phi_{t}^{(1)}]\circ\Delta_{t}%
\circ[\phi_{t}^{(1)}]^{-1},$ and an equivalence $\phi_{t}^{(1)}:H_{t}%
^{(1)}\sim H_{t}$. An easy calculation gives $d_{1}^{(1)}=$ $\mu_{1}%
^{(1)}=\Delta_{1}^{(1)}=0.$ Inductively, given a deformation $H_{t}%
^{(n)}=(H\left[  \left[  t\right]  \right]  ,d_{t}^{(n)},\mu_{t}^{(n)}%
,\Delta_{t}^{(n)})$ with $d_{i}^{(n)}=$ $\mu_{i}^{(n)}=\Delta_{i}^{(n)}=0$ for
$1\leq i\leq n$, and an equivalence $\phi_{t}^{(n)}=\prod_{i=1}^{n}%
(1-t^{i}\phi_{i}):H_{t}^{(n)}\sim H_{t},$ choose $\phi_{n+1}\in\widetilde
{C}^{1}(H;H;3)$ such that $D_{C}^{0,1,1}\left(  \phi_{n+1}\right)  =$
$d_{n+1}^{(n)}+\mu_{n+1}^{(n)}+\Delta_{n+1}^{(n)}$. There is a deformation
$H_{t}^{(n+1)}=(H\left[  \left[  t\right]  \right]  ,d_{t}^{(n+1)},\mu
_{t}^{(n+1)},\Delta_{t}^{(n+1)})$ with $d_{i}^{(n+1)}=$ $\mu_{i}%
^{(n+1)}=\Delta_{i}^{(n+1)}=0$ for $1\leq i\leq n+1$, and an equivalence
$\phi_{t}^{(n+1)}=\prod_{i=1}^{n+1}(1-t^{i}\phi_{i}):H_{t}^{(n+1)}\sim H_{t}$.
Hence there is a sequence $\{\phi_{t}^{(n)}:H_{t}^{(n)}\sim H_{t}\}_{n\geq1}$
that converges $t$-adically to an equivalence $\phi_{t}^{(\infty)}=\prod
_{i=1}^{\infty}(1-t^{i}\phi_{i}):H_{0}\left[  \left[  t\right]  \right]  \sim
H_{t}$ and we conclude that:

\begin{theorem}
If $\widetilde{H}_{d.g.h.a.}^{2}(H;H;3)=0$ then $(H,d,\mu,\Delta)$ is rigid as
a \textit{d.g.h.a.}
\end{theorem}

Furthermore, since the obstructions to extending some finite approximation to
a deformation lie in $\widetilde{H}_{d.g.h.a.}^{3}(H;H;3)$ we have:

\begin{theorem}
If $\widetilde{H}_{d.g.h.a.}^{3}(H;H;3)=0$ then every linear approximation
extends to a deformation$.$
\end{theorem}

\noindent Thus $\widetilde{H}_{d.g.h.a.}^{*}(H;H;3)$ directs the deformation
theory of $H_{0}$ as a \textit{d.g.h.a.} in the direction of the
''infinitesimals'', i.e., $2$-cocycles. In the sequel, to which we refer the
reader, Lazarev and Movshev \cite{Lazarev} apply this particular theory to
analyze the deformations of the de Rham cochains on a Lie group\textit{.}

Finally, we note that the appropriate setting for the deformation theory of
Hopf algebras as quasi-Hopf algebras \cite{Drinfeld1}, \cite{Drinfeld2} is
$\{C^{0,m,n}(H;H),\partial_{C},\delta_{C}\}_{m\geq1;\text{ }n\geq0},$ i.e.,
the subcomplex of cochains in the ''semi-restricted coordinate plane'' $p=0,$
$m\geq1.$\bigskip\ \bigskip\ \ 

I must express my sincere appreciation to Jim Stasheff who proposed this
project and offered many helpful comments and suggestions along the way. Also,
I wish to thank Don Schack for his helpful critique of the manuscript.

\end{document}